\documentclass[10pt,a4paper]{article}
\usepackage[latin1]{inputenc}
\usepackage{amsmath}
\usepackage{amsfonts,amssymb}
\usepackage{epsfig}
\usepackage{psfrag}
\usepackage{mathrsfs}

\def\dsp{\displaystyle}

\def\R{\mathbb{R}}

\def\eps{\varepsilon}

\def\b0{\boldsymbol{0}}
\def\bA{\boldsymbol{A}}
\def\bB{\boldsymbol{B}}
\def\bC{\boldsymbol{C}}
\def\bE{\boldsymbol{E}}
\def\bF{\boldsymbol{F}}
\def\bg{\boldsymbol{g}}
\def\bG{\boldsymbol{G}}
\def\bH{\boldsymbol{H}}
\def\bn{\boldsymbol{n}}

\def\bu{\boldsymbol{u}}
\def\bv{\boldsymbol{v}}
\def\bx{\boldsymbol{x}}
\def\bX{\boldsymbol{X}}
\def\bdelta{\boldsymbol{\delta}}
\def\bsigma{\boldsymbol{\sigma}}
\def\btau{\boldsymbol{\tau}}

\def\Re{\mathfrak{Re}}
\def\We{\mathfrak{We}}

\def\div{\mbox{\,{\textrm{div}}}}

\def\transp{\,^T\!}
\newcommand{\proof}{\noindent{\bf Proof~~}}
\newcommand{\cqfd}{{\nobreak\hfil\penalty50\hskip2em\hbox{}\nobreak\hfil $\square$\qquad\parfillskip=0pt\finalhyphendemerits=0\par\medskip}}

\newtheorem{theorem}{Theorem}[section]
\newtheorem{lemma}{Lemma}[section]
\newtheorem{remark}{Remark}[section]

\newtheorem{proposition}{Proposition}[section]

\begin{document}

\begin{center}
{ \Large \bf Existence results for integral viscoelastic fluids}
\end{center}

\begin{center}
Laurent Chupin
\footnote{
Laboratoire de Mathématiques, UMR 6620 \\
Universit\'e Blaise Pascal, Campus des C\'ezeaux, F-63177 Aubi\`ere Cedex, France.\\
{\it laurent.chupin@math.univ-bpclermont.fr}
}
\end{center}

\begin{abstract}
\noindent We consider the flows of viscoelastic fluid which obey a constitutive law of integral type.
The existence and uniqueness results for solutions of the initial boundary value problem are proved, and the stationary case is studied.
\end{abstract}



\section{Introduction}\label{part:introduction}

The objective of this paper is to provide mathematical results on the integral models for viscoelastic flows: existence solution, uniqueness and stationary study.
\par\noindent
In integral models, stress components $\btau$ are obtained by integrating appropriate functions, representing the amount of deformation, over the strain history of the fluid.
The constitutive law given the stress at time~$t$ and at position~$\bx$ is written as follows:
\begin{equation}\label{eq:2219i}
\btau(t,\bx) = \int_{-\infty}^{t} m(t-T) \, \mathcal S\big(\bF(T,t,\bx)\big)\, \mathrm dT.
\end{equation}
The scalar function $m$ (the memory) and the tensorial function $\mathcal S$ are given by the properties of the fluids studied, whereas the deformation tensor $\bF$ is coupled with the velocity field of the flow.
This flow is itself governed by the Navier-Stokes equations, this constitutes a very strong coupling between the velocity and the stress.
\par\noindent
In some rare cases, it is possible to express the integral models into differential forms (as in the Maxwell models which corresponds to the case where~$\mathcal S$ is linear and where the memory~$m$ exponentially decreases).
For these differential models, there are many mathematical results in the same spirit of those presented here (see for instance~\cite{Chupin2,Fernandez-Guillen-Ortega,Guillope-Saut1,Guillope-Saut3,Lions-Masmoudi-viscoelastique,Molinet}).
\par\noindent
But for really integral models there are far fewer results.
One of the only relevant work on this type of model is that of M.~Renardy~\cite{Renardy}.
In its paper, M.~Renardy prove an existence and uniqueness result for a K-BKZ fluid using Kato's theory of quasilinear hyperbolic equations.
Its elegant approach differs substantially from the approach used here, and does not seem easily adaptable to more general laws.
For instance, in~\cite{Renardy} the author assumes that~$m$ is not singular at~$0$ contrary to what is predicted by some molecular models like the Doi-Edwards model.
It is important to note that while the theoretical results are very few, many authors have studied the numerical simulation of flow with an integral law of type~\eqref{eq:2219i}.
The review article~\cite{Keunings} and references cited therein, provide a good overview of the state of the art regarding the various methods.
\par\noindent
The main reason for this lack of theoretical results is probably due to the nature of the equations:
\par
$\bullet$ To evaluate the stress in an integral model, we must know all the previous configurations.
This difficulty is overcome by introducing an additional time variable corresponding to the age.
It is then necessary to manage two different times, and in particular use a Gronwall type lemma in two variables in order to obtain fine estimates with respect to these variables.
\par
$\bullet$ The usual integral models are strongly nonlinear (in the linear case we find the well-known Maxwell models).
The possibility to circumvent this difficulty is to work with solutions regular enough, more precisely in a Banach algebra like the Sobolev spaces~$W^{1,p}$ for $p$ large enough, typically $p$ greater than the dimension of the physical fluid domain.\\[0.3cm]
{\bf Organization of the paper --}
Section~\ref{part:model} is devoted to the presentation of the model.
The dimensionless form of equations, as well as many classic examples are given.
The main results are stated in the Section~\ref{part:results} whereas the proofs are given in the next sections.
Section~\ref{part:local} is entirely devoted to the proof of the first theorem regarding the local (in time) well posseness.
The three next sections~(\ref{part:unic}, \ref{part:global} and~\ref{part:station}) concern the proof of the uniqueness result, the global existence with small data, and the case of stationary solutions respectively.
The conclusion of this paper (section~\ref{part:conclusion}) contains many remarks and open questions.
Finally, some notions on tensors, and a technical Gronwall type lemma has been postponed to Appendices~\ref{part:appendix0} and~\ref{part:appendix1}.

\section{Governing equations}\label{part:model}

\subsection{Conservation principles}\label{part:model-1}

The fluid flows is modeled using the equation of conservation of the linear momentum and the equation of the conservation of mass, which read in the incompressible and isothermal case as follows:
\begin{equation}\label{eq:0001}
\left\{
\begin{aligned}
& \rho \, ( \partial_t \bu + \bu\cdot \nabla \bu ) + \nabla p = \div \, \bsigma + \boldsymbol{f}, \\
& \div \, \bu = 0,
\end{aligned}
\right.
\end{equation}
The real~$\rho$ is a constant density of mass and the vector~$\boldsymbol{f}$ corresponds to some external body forces.
This system is closed using a constitutive equation connecting the stress and the deformation $D\bu = \frac{1}{2}(\nabla \bu + \transp{(\nabla \bu)})$.
\par\noindent
For a so-called Newtonian viscous fluid, the relationship is linear: $\btau_s = 2\eta_s D\bu$.
The real $\eta_s>0$ is named the solvent viscosity and the contribution $\div\,\btau_s$ in the momentum equation gives the usual diffusive term~$\eta_s \Delta \bu$.
To taking account some elasticity aspect, appearing for instance in polymer solution, we add to the viscous contribution~$\btau_s$ an elastic one:
\begin{equation}\label{eq:0002}
\bsigma = 2\eta_s D\bu + \btau.
\end{equation}
The role of this additional contribution~$\btau$ is to take into account the past history of the fluid.
The most natural way to do this is to introduce the integral models.

\subsection{Integral models}

Very generally, the elastic contribution~$\btau(t,\bx)$ at time $t$ and at spatial position~$\bx$ is written 
\begin{equation}\label{eq:2129}
\btau(t,\bx) = \underset{T<t}{\mathfrak F} (\bF(T,t,\bx)),
\end{equation}
where $\mathfrak F$ is a functional to clarify, which depends on the deformation gradient~$\bF(T,t,\cdot)$ from a times $T$ to a next time~$t$.
\par\noindent
More precisely, the deformation gradient~$\bF(T,t,\cdot)$ measures stretch and rotation.
It is defined as follows:
for two times $T\leq t$ given, we first introduce the notation $\bx(T,t,\bX)$ which corresponds to the position at time~$t$ of the fluid particle which was at the position~$\bX$ at time~$T$.
The dynamics of any mechanical problem with a velocity field $\bu(t,\bx)$ can be described by this flow map $\bx(T,t,\bX)$ which is a time dependent family of orientation preserving diffeomorphisms:
\begin{equation}\label{eq:1616}
\left\{
\begin{aligned}
& \partial_t \bx(T,t,\bX) = \bu(t,\bx(T,t,\bX)), \\
& \bx(T,T,\bX) = \bX.
\end{aligned}
\right.
\end{equation}
The deformation gradient $\widetilde \bF(T,t,\bX)$ is used to describe the changing of any configuration, amplification or pattern during the dynamical process.
It is defined by
\begin{equation}\label{eq:1612}
\widetilde \bF(T,t,\bX) = \frac{\partial \bx}{\partial \bX}(T,t,\bX).
\end{equation}
Finally, the deformation gradient $\bF(T,t,\bx)$ will be defined as the corresponding in the Eulerian coordinates:
\begin{equation}\label{eq:1613}
\bF(T,t,\bx(T,t,\bX)) = \widetilde \bF(T,t,\bX).
\end{equation}
The integral models we study in this article correspond to the particular case of equation~\eqref{eq:2129}. They are written
\begin{equation}\label{eq:2219}
\btau(t,\bx) = \int_{-\infty}^{t} m(t-T) \, \mathcal S\big(\bF(T,t,\bx)\big)\, \mathrm dT,
\end{equation}
where $m$ is called a memory function and $\mathcal S$ is a model-dependent strain measure. It is a tensorial function (its arguments and its images are $2$-tensors).
\begin{remark}\label{rem:2010}
Due to many physical principle, the functions $m$ and $\mathcal S$ must satisfy some assumptions.
\par\noindent
$\bullet$ For instance, the principle of frame indifference implies that the stress tensor depends on the relative deformation gradient~$\bF$ only through the relative Finger tensor~$\transp{\bF} \cdot \bF$ (or its inverse, the Cauchy-Green tensor), see the examples given in the subsection~\ref{part:example}.
\par\noindent
$\bullet$ In the same way (see also the subsection~\ref{part:example}), the principle of fading memory implies that $m$ must be a positive function which decreases to~$0$.
\end{remark}

\subsection{A closed system}\label{part:model-3}

\subsubsection{PDE for the deformation gradient}

We derivate the relation~\eqref{eq:1613} with respect to the time~$t$.
The chain rule together with the relation~\eqref{eq:1616} yields the following equation
\begin{equation}\label{eq:1624}
\begin{aligned}
\partial_t \widetilde \bF(T,t,\bX) 
& = \partial_t \bF (T,t,\bx) + \partial_t \bx(T,t,\bX) \cdot \partial_{\bx} \bF (T,t,\bx) \\
& = \partial_t \bF (T,t,\bx) + \bu(t,\bx) \cdot \partial_{\bx} \bF (T,t,\bx).
\end{aligned}
\end{equation}
But using the relation~\eqref{eq:1612} together with the chain rule and the relation~\eqref{eq:1616} again, we obtain
\begin{equation}\label{eq:1625}
\begin{aligned}
\partial_t \widetilde \bF(T,t,\bX) 
& = \partial_{\bX} (\partial_t \bx(T,t,\bX)) \\
& = \partial_{\bX} ( \bu(t,\bx(T,t,\bX)) ) \\
& = \partial_{\bX} \bx(T,t,\bX) \cdot \partial_{\bx}\bu(t,\bx) \\
& = \bF(T,t,\bx) \cdot \partial_{\bx}\bu(t,\bx).
\end{aligned}
\end{equation}
The equations~\eqref{eq:1624} and~\eqref{eq:1625} show that we have the following relation coupling the velocity field~$\bu$ and the deformation gradient~$\bF$:
\begin{equation}\label{eq:1650}
\partial_t \bF + \bu \cdot \nabla \bF = \bF \cdot \nabla \bu.
\end{equation}

\subsubsection{A new time variable to take into account the past}

Note that in the previous subsection, the time $T$ can be view as a parameter.
In fact, it is only used in the law~\eqref{eq:2129}, or in the law~\eqref{eq:2219} for the integral form, as a marker of past events.
In the sequel, it is interesting to select as independent variable the age $s=t-T$, which is measured relative to the current time~$t$.
We now introduce $\bG(s,t,\bx)=\bF(t-s,t,\bx)$.
Clearly, we have the following relation instead of the relation~\eqref{eq:1650}:
\begin{equation}\label{eq:1651}
\partial_t \bG + \partial_s \bG + \bu \cdot \nabla \bG = \bG \cdot \nabla \bu,
\end{equation}
where naturally the velocity $\bu$ only depends on $(t,\bx)$, and is independent of this new variable~$s$.
Moreover, in term of variables $(s,t)$, the relation~\eqref{eq:2219} given the stress tensor reads
\begin{equation}\label{eq:2220}
\btau(t,\bx) = \int_0^{+\infty} m(s) \, \mathcal S\big(\bG(s,t,\bx)\big)\, \mathrm ds.
\end{equation}
\noindent
{\bf Initial and past conditions --}
The equation~\eqref{eq:1651} requires initial conditions, that is to say that we must give the value $\bG(s,0,\bx)$ at time $t=0$ for all age~$s\geq 0$ and for all point~$\bx\in \Omega$, and the value $\bG(0,t,\bx)$ at any time $t\geq 0$ and for all point~$\bx\in \Omega$.
Note that for $T$ fixed, the deformation field $\bF(T,t,\bx)$ can be thought of as having been created at time $t = T$ with the natural initial condition $\bF(T,T,\bx)=\bdelta$.
In term of the new variables $(s,t)$ this relation reads
\begin{equation}
\bG(0,t,\bx)=\bdelta.
\end{equation}
The condition about $\bG(s,0,\bx)$ correspond to the description of the deformation gradient before the time $t=0$.
In practice this value is unknown but it is physically reasonable to assume that for an age old enough, the fluid was quiescent.
In practice we will only assume that there exists a given function~$\bG_{\mathrm{old}}$ such that
\begin{equation}
\bG(s,0,\bx)=\bG_{\mathrm{old}}(s,\bx).
\end{equation}

\subsubsection{Non-dimensional final model}

The resulting system using the equations~\eqref{eq:0001}, \eqref{eq:0002}, \eqref{eq:1651} and~\eqref{eq:2220} is then written
\begin{equation}\label{eq:0003}
\left\{
\begin{aligned}
& \rho \, ( \partial_t \bu + \bu\cdot \nabla \bu ) + \nabla p - \eta_s \Delta \bu = \div \, \btau + \boldsymbol{f}, \\
& \div \, \bu = 0, \\
& \btau(t,\bx) = \int_0^{+\infty} m(s) \, \mathcal S\big(\bG(s,t,\bx)\big)\, \mathrm ds, \\
& \partial_t \bG + \partial_s \bG + \bu \cdot \nabla \bG = \bG \cdot \nabla \bu.
\end{aligned}
\right.
\end{equation}
This system can be adimensionalized in the usual way.
We introduce  the characteristic values $U$ and $L$ for the velocity and the length. The current time is then of order of $L/U$ and it is natural to introduce another characteristic time~$\lambda$ for the age variable~$s$.
The characteristic viscosity of the fluid takes into account the viscosity $\eta_s$ of the solvent, but also the viscosity~$\eta_e$ of the elastic part (the polymer): we note $\eta=\eta_s+\eta_e$.
More precisely, we introduce the following dimensionless variable, quoted by a star:
\begin{equation}
\begin{aligned}
& \bx^\star = \frac{\bx}{L},
\quad
\bu^\star = \frac{\bu}{U},
\quad
t^\star = \frac{t}{L/U},
\quad
s^\star=\frac{s}{\lambda},
\quad
p^\star = \frac{p}{\eta\, U/L},
\quad
\boldsymbol f^\star = \frac{\boldsymbol f}{\eta \, U/L^2},\\
& \btau^\star = \frac{\btau}{\eta \, U/L},
\quad
\bG^\star = \bG,
\quad
\mathcal S^\star(\bG^\star) = \frac{\mathcal S(\bG)}{\eta_e/\lambda},
\quad
m^\star(s^\star) = \frac{m(s)}{1/\lambda}.
\end{aligned}
\end{equation}
The system~\eqref{eq:0003} is then written in dimensionless form as follows (where we drop the star for sake of simplicity):
\begin{equation}\label{eq:1321}
\left\{
\begin{aligned}
& \Re ( \partial_t \bu + \bu\cdot \nabla \bu ) + \nabla p - (1-\omega) \Delta \bu = \div \, \btau + \boldsymbol f, \\
& \div \, \bu = 0, \\
& \btau(t,\bx) = \frac{\omega}{\We} \int_0^{+\infty} m(s) \, \mathcal S\big(\bG(s,t,\bx)\big)\, \mathrm ds,\\
& \partial_t \bG + \frac{1}{\We} \partial_s \bG + \bu \cdot \nabla \bG = \bG \cdot \nabla \bu,
\end{aligned}
\right.
\end{equation}
where we introduced the three non-dimensional numbers which characterize the flow:
\begin{itemize}
\item[$\bullet$] The Reynolds number $\Re = \frac{\rho\, U L}{\eta}$ which corresponds to the ratio between inertial and viscous forces acting on the fluid;
\item[$\bullet$] The Weissenberg number $\We = \frac{\lambda\, U}{L}$ which is the ratio between the time of the relaxation of the fluid and the time of the experiment;
\item[$\bullet$] The retardation parameter $\omega = \frac{\eta_e}{\eta} \in [0,1]$ which balances the purely viscous effects ($\omega=0$) and the purely elastic effects ($\omega=1$).
\end{itemize}
The system~\eqref{eq:1321} is closed with the following initial and boundary conditions:
\begin{equation}\label{eq:1322}
\bu\big|_{t=0} = \bu_0, \qquad \bu\big|_{\partial \Omega} = \b0, \qquad
\bG\big|_{t=0} = \bG_{\mathrm{old}}, \qquad \bG\big|_{s=0} = \bdelta.
\end{equation}
\begin{remark}\label{rem:2011}
As we precised in the Remark~\ref{rem:2010}, the stress tensor~$\btau$ depends on the deformation tensor $\bG$ {\it via} the right relative Finger tensor~$\bB = \transp{\bG} \cdot \bG$ or {\it via} the Green-Cauchy tensor~$\bC = \bB^{-1}$.
Using the last equation of~\eqref{eq:1321}, we note that the tensor~$\bB$ satisfies
\begin{equation}\label{eq:1652a}
\partial_t \bB + \frac{1}{\We} \partial_s \bB + \bu \cdot \nabla \bB = \bB \cdot \nabla \bu + \transp{(\nabla \bu)}\cdot \bB,
\end{equation}
whereas the tensor~$\bC$ satisfies
\begin{equation}\label{eq:1652b}
\partial_t \bC + \frac{1}{\We} \partial_s \bC + \bu \cdot \nabla \bC = - \bC \cdot \transp{(\nabla \bu)} - \nabla \bu \cdot \bC.
\end{equation}
\end{remark}

\subsection{Examples of integral models}\label{part:example}

In this section, we present some classical integral laws of kind~\eqref{eq:2220} to model the viscoelasticity.
These law are defined by the memory function~$m$ and by the strain measure~$\mathcal S$.

\subsubsection{Memory function~$m$}\label{part:example-0}

Usually, in the viscoelastic formulas, a relaxation function~$g$ is introduced to describe the return into equilibrium.
It physically corresponds to the response of the stress to a shear jump.
In accordance with thermodynamics through what is called the principle of fading memory - see~\cite{Coleman} - this relaxation function is constrained in that $g(0) = 1$ and $\lim_{+\infty} g = 0$.
To so called memory function~$m$ used in the integral models~\eqref{eq:2220} corresponds to $m = -g'$.
Consequently, the function~$m$ is decreasing, positive and satisfies $\int_0^\infty m(s)\, \mathrm ds = 1$.
\par\noindent
In many cases experimentally observed relaxation functions exhibit a stretched exponential decay $\mathrm e^{-(s/\lambda)}$ where $\lambda>0$ is a relaxation time.
We could as well have considered the case of several relaxation times, that is a memory function like
\begin{equation}\label{ex:m1}
m(s) = \sum_{k=1}^N \frac{\eta_k}{\lambda_k^2} \mathrm e^{-s/\lambda_k}.
\end{equation}
From a mathematical point of view, it will be equivalent to consider $m(s)=\mathrm e^{-s}$ (on the dimensionless form).
\par\noindent
This expression for the memory function can be generalized.
For instance, in the Doi-Edwards model - see~\cite{Doi}, the memory function is given by
\begin{equation}\label{ex:m2}
m(s)=\frac{8}{\pi^2 \, \lambda} \sum_{k=0}^{+\infty} \mathrm e^{-(2k+1)^2s/\lambda}.
\end{equation}
In practice, the difference between this model~\eqref{ex:m2} and the model~\eqref{ex:m1} containing a finite number of relaxation times is really important.
In fact, in the case of the model~\eqref{ex:m2}, the function~$m$ has a singularity in~$0$.
This singularity can bring additional difficulties (eg, such a case is not treated in the article~\cite{Renardy} of M.~Renardy).
The function~$m$ remains integrable, which is the key assumption for the present results.\\[0.3cm]
Even if the exponential case is usually used, many other possible choices for the memory function~$m$ are possibles - see~\cite{Freed-Diethelm}.
The algebraic pattern $g(s)=(s/\lambda)^{-\beta}$ with $0<\beta<1$ is observed in the stress relaxation of viscoelastic materials such as critical gels~\cite{Chambon,Schiessel}, in the charge carrier transport in amorphous semiconductors~\cite{Scher}, in dielectric relaxation~\cite{Jonscher} or in the attenuation of seismic waves~\cite{Kjartansson}.
That corresponds to the following memory functions
\begin{equation}
m(s) = \sum_{k=1}^N \frac{\eta_k\, \beta_k}{\lambda_k} \Big( \frac{s}{\lambda_k} \Big)^{-(\beta_k+1)}.
\end{equation}

\subsubsection{Strain measure~$\mathcal S$}\label{part:example-1}

{\bf Maxwell models --}
The more simple case corresponds to the choice $\mathcal S(\bG)=\bB-\bdelta$ where $\bB=\transp{\bG}\cdot \bG$, and where the memory (dimensionless) function~$m$ is given by $m(s)=\mathrm e^{-s}$.
The stress tensor $\btau$ is then given by
\begin{equation}
\btau(t,\bx) = \frac{\omega}{\We} \int_0^{+\infty} \mathrm e^{-s} \big( \bB(s,t,\bx)-\bdelta \big) \, \mathrm ds.
\end{equation}
This expression is simple enough to deduce a PDE for the stress tensor $\btau$ from the PDE for the deformation tensor $\bG$.
In fact we use the equation~\eqref{eq:1652a} satisfied by the Finger tensor:
\begin{equation}\label{eq:1654}
\partial_t (\bB-\bdelta) + \frac{1}{\We} \partial_s (\bB-\bdelta) + \bu \cdot \nabla (\bB-\bdelta) = (\bB-\bdelta) \cdot \nabla \bu + \transp{(\nabla \bu)}\cdot (\bB-\bdelta) + 2 D\bu.
\end{equation}
We next multiply this equation by $\frac{\omega}{\We} m(s)$ and integrate for $s\in (0,+\infty)$.
Taking into account the initial condition $\bB|_{s=0} = \bdelta$, we obtain the Upper Convected Maxwell (UCM) model:
\begin{equation}
\We \Big( \partial_t \btau + \bu \cdot \nabla \btau - \transp{(\nabla \bu)} \cdot \btau - \btau \cdot \nabla \bu \Big) + \btau = 2 \omega D\bu.
\end{equation}
Another classical case is to the Lower Convected Maxwell (LCM) model.
It corresponds to $\mathcal S(\bG) = \bdelta - \bC$ where $\bC=\bB^{-1}$, and to an exponential memory function $m(s)=\mathrm e^{-s}$.
Using the equation~\eqref{eq:1652b} we obtain, like to get the UCM model:
\begin{equation}
\We \Big( \partial_t \btau + \bu \cdot \nabla \btau + \btau \cdot \transp{(\nabla \bu)} + \nabla \bu \cdot \btau \Big) + \btau = 2 \omega D\bu.
\end{equation}

\begin{remark}
In fact there exists a continuum of such model (called Oldroyd models, corresponding to a balance between upper-convected and lower-convected) but we do not know if these models derive from integral models.
\end{remark}

\noindent
{\bf K-BKZ models --}
Among the most relevant non-linear cases, the most popular integral models for a viscoelastic flow are the K-BKZ models introduce by B. Bernstein, E. A. Kearsley and L. J. Zapas~\cite{KBKZ1, KBKZ2} and A. Kaye~\cite{Kaye}.
For such models, $\mathcal S$ takes the following form:
\begin{equation}
\mathcal S(\bG) = \phi_1(I_1,I_2) (\bB - \bdelta) + \phi_2(I_1,I_2) (\bdelta - \bC),
\end{equation}
where $\phi_1$ and $\phi_2$ are two scalar functions of the strain invariants~$I_1=\mathrm{Tr}(\bB)$ and~$I_2=\mathrm{Tr}(\bB^{-1})$ (see the Appendix~\ref{part:appendix0} for a discussion on these invariants).
Clearly, the two Maxwell models presented in the previous paragraph are particular K-BKZ models, for $(\phi_1,\phi_2)=(1,0)$ and $(\phi_1,\phi_2)=(0,1)$.

\begin{remark}
Add a diagonal tensor to $\mathcal S(\bG)$, that is consider $\mathcal S(\bG)+\phi \, \bdelta$ instead of $\mathcal S(\bG)$, do not change the mathematical structure of the system.
In fact $\btau$ becomes $\btau+\frac{\omega}{\We} \phi \, \bdelta$ and the additional contribution can be consider as a pressure contribution.
\end{remark}

\noindent
Following the last remark, we can see the PSM models presented by A. C.~Papanastasiou, L.~Scriven and C.~Macosko in~\cite{Papanastasiou} as K-BKZ models:
\begin{equation}\label{ex:0001}
\mathcal S(\bG) = h(I_1,I_2) \bB
\qquad \text{with} \quad h(I_1,I_2)  = \frac{\alpha}{\alpha + \beta I_1 + (1-\beta) I_2 - 3}.
\end{equation}
In these models, the parameters $\alpha>0$ and~$0\leq \beta \leq 1$ are obtained from the rheological fluid properties.
In the same way, Wagner~\cite{Wagner0,Wagner1} proposes the following law:
\begin{equation}\label{ex:0002}
\mathcal S(\bG) = h(I_1,I_2) \bB
\qquad \text{with} \quad h(I_1,I_2)  = \mathrm{exp}(-\alpha\sqrt{ \beta I_1 + (1-\beta) I_2 - 3}).
\end{equation}

\noindent
{\bf Doi-Edwards model --}
The Doi-Edwards model is a molecular model where the motion of the polymers is described by reptation in a tube, more precisely it corresponds to the simplest tube model of entangled linear polymers.
The memory function~$m$ associated to such model is given by the relation~\eqref{ex:m2} whereas the strain measure is obtained as an average with respect to the orientation of tube segments.
Works of P.-K. Currie show that we can approach this model using the following strain function (named the Currie approximation, see~\cite{Currie}):
\begin{equation}\label{ex:0003}
\mathcal S(\bG) = \frac{4}{3(J-1)}\bB - \frac{4}{3(J-1)\sqrt{I_2+3.25}}\bC
\qquad \text{with} \quad J=I_1+2\sqrt{I_2+3.25}.
\end{equation}

\section{Main results}\label{part:results}

\subsection{Mathematical framework}

In the sequel, the fluid domain $\Omega\subset \R^d$, $d\geq 2$, is a bounded connected open set with a Lipschitz continuous boundary $\partial \Omega$.
We use the following standard notations:
\begin{itemize}
\item For all real $s\geq 0$ and all integer $p\geq 1$, the set $W^{s,p}(\Omega)$ corresponds to the usual Sobolev spaces.
We classically note $L^p(\Omega)=W^{0,p}(\Omega)$ and $H^s(\Omega)=W^{s,2}(\Omega)$.
\end{itemize}
We will frequently use functions with values in $\R^d$ or in the space~$\mathcal L(\R^d)$ of real $d\times d$ matrices.
In all cases, the notations will be abbreviated.
For instance, the space $(W^{1,p}(\Omega))^3$ will be denoted $W^{1,p}(\Omega)$.
Moreover, all the norms will be denoted by index, for instance like $\|\bu\|_{W^{1,p}(\Omega)}$.
\begin{itemize}
\item Since we are interested in the incompressible flows, we introduce
$
H_p(\Omega) = \{\bv \in L^p(\Omega) \, ; \, \div\, \bv=0,\, \bv\cdot \bn=0 \text{ on }\partial\Omega\},
$
where $\bn$ is the unitary vector normal to $\partial \Omega$, oriented towards the exterior of~$\Omega$.
Moreover, we note $V(\Omega)=H_2(\Omega)\cap H^1_0(\Omega)$ and~$V'(\Omega)$ its dual.
\item The Stokes operator $A_p$ in $H_p(\Omega)$ is introduced, with domain $D(A_p)(\Omega) = W^{2,p}(\Omega)\cap W^{1,p}_0(\Omega) \cap H_p(\Omega)$, whereas we note $D^r_p(\Omega) = \{\bv \in H_p \, ; \, \|\bv\|_{L^p(\Omega)} + \big( \int_0^{+\infty} \|A_p \mathrm e^{-tA_p}\bv \|_{L^p(\Omega)}^{r} \, \mathrm dt \big)^{1/r} < +\infty \}$.
\item The notation of kind $L^r(0,T;D(A_p))$ denotes the space of $r$-integrable functions on $(0,T)$, $T>0$, with values in~$D(A_p)$.
Similarly, expressions like $g\in L^\infty(\R^+;L^r(0,T;L^p(\Omega)))$ means that
$$
\sup_{s\in \R^+} \Big( \int_0^T \|g(s,t,\cdot)\|_{L^p(\Omega)}^r \, \mathrm dt \Big)^{\frac{1}{r}} < +\infty.
$$
\end{itemize}
{\bf Assumptions --}
About the model~\eqref{eq:1321} itself, it uses the two given functions~$m$ and~$\mathcal S$.
From a mathematical point of view we assume very general assumptions (satisfied by all the physical models introduced earlier):
\begin{itemize}
\item[(A1)] $m:s\in \R^+\longmapsto m(s)\in \R$ is measurable, positive and satisfies $\dsp \int_0^{+\infty} m(s) \, \mathrm ds = 1$;
\item[(A2)] $\mathcal S:\bG\in \mathcal L(\R^d) \longmapsto \mathcal S(\bG) \in \mathcal L(\R^d)$ is of class $\mathcal C^1$.
\end{itemize}
Note that the notion of derivative for the $2$-tensorial application~$\mathcal S$ will be precised in the Appendix~\ref{part:appendix0}.
Moreover, if we wanted to be more precise, the assumption (A2) is written rather ``the function~$\mathcal S$ is of class~$\mathcal C^1$ on a subset of $\mathcal L(\R^d)$ taking into account the fact that $\det \bG = 1$'' - see Appendix~\ref{part:appendix0} again.
Throughout the remainder of this paper, these two hypotheses (A1) and (A2) will be assumed satisfied.

\subsection{Statements of main results}

The first result concerns an existence result for strong solutions.
It is obviously local with respect to time (as for the results on the Navier-Stokes equations):
\begin{theorem}[local existence]\label{th:local}
Let $T>0$, $r\in ]1,+\infty[$ and $p\in]d,+\infty[$.
\par\noindent
If $\bu_0 \in D^r_p(\Omega)$, $\bG_{\mathrm{old}}\in L^\infty(\R^+;W^{1,p}(\Omega))$, $\partial_s \bG_{\mathrm{old}}\in L^r(\R^+;L^p(\Omega))$ and $\boldsymbol{f} \in L^r(0,T;L^p(\Omega))$ then there exists $T_\star \in ]0,T]$ and a strong solution $(\bu,p,\btau,\bG)$ to the system~\eqref{eq:1321} in $[0,T_\star]$, which satisfies the initial/boundary conditions~\eqref{eq:1322}.
Moreover we have
\begin{equation}
\begin{array}{ll}
 \bu \in L^r(0,T_\star;W^{2,p}(\Omega)),
& \partial_t \bu \in L^r(0,T_\star;L^p(\Omega)),\\
 \btau \in L^\infty(0,T_\star;W^{1,p}(\Omega)),
& \partial_t \btau \in L^r(0,T_\star;L^p(\Omega)),\\
 \bG \in L^\infty(\R^+ \! \times \! (0,T_\star);W^{1,p}(\Omega)), \qquad
& \partial_s \bG, \, \partial_t \bG \in L^\infty(\R^+;L^r(0,T_\star;L^p(\Omega))).
\end{array}
\end{equation}
\end{theorem}
We will show that the solution obtained in the theorem~\ref{th:local} is the only one in the class of regular solution.
Precisely, the result reads as follow
\begin{theorem}[uniqueness]\label{th:unic}
Let $T>0$.
\par\noindent
If $\bu_0\in H$, $\bG_{\mathrm{old}}\in L^\infty(\R^+;L^2(\Omega))$, $\boldsymbol{f}\in L^1(0,T;V'(\Omega))$ and if \eqref{eq:1321}-\eqref{eq:1322} possess two weak solutions $(\bu_1,p_1,\btau_1,\bG_1)$ and $(\bu_2,p_2,\btau_2,\bG_2)$ in the usual sense, with for $i\in \{1,2\}$,
\begin{equation}
\begin{aligned}
& \bu_i \in L^\infty(0,T;L^2(\Omega)) \cap L^2(0,T;H^1(\Omega)) \cap L^1(0,T;W^{1,\infty}(\Omega)) \\
& \bG_i \in L^\infty(\R^+ \! \times \! (0,T);W^{1,d}(\Omega))),
\end{aligned}
\end{equation}
then they coincide ($p_1$ and $p_2$ coincide up to an additive function only depending on~$t$).
\end{theorem}
If the data are small, it is possible to show that the unique local solution obtained by the theorem~\ref{th:local} and~\ref{th:unic} is defined for all time (up to an assumption on the relaxation parameter~$\omega$):
\begin{theorem}[global existence with small data]\label{th:global}
Let $r\in ]1,+\infty[$ and $p\in]d,+\infty[$.
\par\noindent
For each $T>0$, with the same assumptions as in the theorem~\ref{th:local}, there exists $\omega_T\in (0,1)$ such that if $0\leq\omega<\omega_T$ and if the data $\bu_0$ and $\boldsymbol{f}$ have sufficiently small norms in their respective spaces, then there exists a unique strong solution $(\bu,p,\btau,\bG)$ to system~\eqref{eq:1321}-\eqref{eq:1322} in $[0,T]$ which belongs in the same spaces that the local solution obtained in the theorem~\ref{th:local}.
\end{theorem}
The last theorem which is proved in this paper concerns the stationary solutions, that is solutions which do not depend on time.
Since the mathematical problem is fundamentally different, the proof that we propose (see Section~\ref{part:station}) requires some additional assumptions on the memory function~$m$ and on the function~$\mathcal S$.
\begin{itemize}
\item[(A3)] There exists $\alpha>0$ and $c\geq 0$ such that $m(s)\leq c \, \mathrm e^{-\alpha s}$ for all $s\in \R^+$.
\item[(A4)] The function $\mathcal S$ and its derivative have polynomial growth: there exists $(a,b,c)\in \R^3$ such that for all $\bG\in \mathcal L(\R^d)$ we have $|\mathcal S(\bG)| \leq c |\bG|^a$ and $|\mathcal S'(\bG)| \leq c |\bG|^b$.
\end{itemize}
Note that the assumption (A3) implies that the memory~$m$ decays sufficiently fast (exponentially) at infinity, and prohibits a singularity at~$0$.
The assumption~(A4) is generally satisfied by classical law (see the Appendix~\ref{part:appendix0} in which we show that the assumption (A4) is satisfied with $a=0$ and $b=-1$ for a PSM law).
\par\noindent
Under these additional assumptions (A3) and (A4), we have:
%
\begin{theorem}[stationary solutions]\label{th:station}
Let $p\in]d,+\infty[$.
\par\noindent
If $\omega$ is small and if $\boldsymbol f \in L^p(\Omega)$ have a sufficiently small norm then there exists exactly one small strong solution $(\bu,p,\btau,\bG)$ to the stationary problem associated to~\eqref{eq:1321}-\eqref{eq:1322} with
\begin{equation}
 \bu \in W^{2,p}(\Omega),
\quad \btau \in W^{1,p}(\Omega),
\quad \bG \in L^\infty(\R^+(\mathrm d\mu);W^{1,p}(\Omega)),
\quad \partial_s\bG \in L^\infty(\R^+(\mathrm d\mu);L^p(\Omega)),
\end{equation}
where $\mathrm d\mu=\mathrm e^{-\beta s}\mathrm ds$, $\beta$ being a positive number (depending on~$p$, $d$, $\omega$ and~$\We$).
\end{theorem}

\begin{remark}\label{rem:theorem}
Here are some remarks concerning the four above theorems:
\par\noindent
- We remark that the theorem~\ref{th:global} do not contain smallness condition on the deformation tensor~$\bG_{\mathrm{old}}$.
At rest this tensor is not zero but is equal to the identity tensor.
An assumption of smallness should eventually be introduced on the quantity~$\bG_{\mathrm{old}}-\bdelta$.
In fact it is implicit in the smallness assumption on the parameter~$\omega$.
In the same way, the term ``small solution'' used in the theorem~\ref{th:station} does not concern the tensor~$\bG$.
\par\noindent
- The theorem~\ref{th:global} can be viewed as a result of stability for the solution $\bu=\b0$, $\btau=\b0$.
\par\noindent
- The assumptions (A3) and (A4) are not optimal (a more precise formulation will be complex).
For instance, following the proof of the theorem~\ref{th:station} (see section~\ref{part:station}) we can see that if $a=b+1=0$ then the assumption~(A3) is unnecessary.
\par\noindent
- An example of value for $a$ and $b$ is given in the Appendix~\ref{part:appendix0} for the PSM model (more precisely for the model~\eqref{ex:0001} with $\alpha=4$ and $\beta=1$).
In this case we have $a=b+1=0$.
\end{remark}

\section{Proof for the local existence}\label{part:local}

This section is devoted to the proof of the local existence theorem~\ref{th:local}.
The main ideas to prove the theorem~\ref{th:local} are based on works of J.C.~Saut and C.~Guillop\'e~\cite{Guillope-Saut-CRAS, Guillope-Saut3}.
Roughly speaking, we rewrite the equations~\eqref{eq:1321} as a fixed point equation and applying the Shauder's theorem.
This principle was taken up by E.~Fernandez-Cara, F.~Guillen and R.R.~Ortega~\cite{Fernandez-Guillen-Ortega-CRAS, Fernandez-Guillen-Ortega} in the context of the functional spaces $L^r\!-\!L^p$.
That is this choice which is presented in the present paper.
\par\noindent
We then first analyze three independent problems.
A linear Stokes system with given source term and initial value:
\begin{equation}\label{eq:1337}
\left\{
\begin{aligned}
& \Re \, \partial_t \bu + \nabla p - (1-\omega) \Delta \bu = \overline \bg, \\
& \div\, \bu = 0,\\
& \bu|_{\partial \Omega} = \b0, \quad \bu|_{t=0} = \bu_0;
\end{aligned}
\right.
\end{equation}
A following tensor problem given the deformation gradient as a function of a given velocity field~$\overline \bu$:
\begin{equation}\label{eq:1338}
\left\{
\begin{aligned}
& \partial_t \bG + \frac{1}{\We} \partial_s \bG + \overline \bu \cdot \nabla \bG = \bG \cdot \nabla \overline \bu, \\
& \bG|_{s=0} = \bdelta,\qquad \bG|_{t=0} = \bG_{\mathrm{old}}; 
\end{aligned}
\right.
\end{equation}
And the constitutive integral law given the stress tensor~$\btau$ with respect to a deformation gradient~$\overline \bG$:
\begin{equation}\label{eq:1339}
\btau(t,\bx) = \frac{\omega}{\We} \int_0^{+\infty} m(s) \, \mathcal S\big(\overline \bG(s,t,\bx)\big)\, \mathrm ds.
\end{equation}

\subsection{Estimates for the velocity~$\bu$ solution of a Stokes problem}

The results for the Stokes system~\eqref{eq:1337} are very usual.
In this subsection we only recall, without proof (we can found a proof in~\cite{Giga}), a well known results for the time dependent Stokes problem:
\begin{lemma}\label{lem:stokes}
Let $T>0$, $r\in ]1,+\infty[$ and $p\in]1,+\infty[$.
\par\noindent
If $\bu_0 \in D^r_p(\Omega)$ and $\overline \bg \in L^r(0,T;H_p)$ then there exists a unique solution $\bu\in L^r(0,T;D(A_p))$ such that $\partial_t\bu \in L^r(0,T;H_p)$ to equations~\eqref{eq:1337}.
Moreover this solution satisfies
\begin{equation}
 \|\bu\|_{L^r(0,T;W^{2,p}(\Omega))} + \|\partial_t \bu\|_{L^r(0,T;L^p(\Omega))} \leq \frac{C_1}{1-\omega} \big( \Re \, \|\bu_0\|_{W^{2,p}(\Omega)} + \|\overline \bg\|_{L^r(0,T;L^p(\Omega))} \big),
\end{equation}
where the constant $C_1$ only depends on $\Omega$, $r$ and $p$.
\end{lemma}

\subsection{Estimates for the deformation gradient~$\bG$}

This subsection is devoted to the proof of the following lemma, which gives estimates for the solution to the system~\eqref{eq:1338}:
\begin{lemma}\label{lem:1117}
Let $T>0$, $r\in ]1,+\infty[$ and $p\in]d,+\infty[$.
\par\noindent
If $\bG_{\mathrm{old}}\in L^\infty(\R^+;W^{1,p}(\Omega))$, $\partial_s\bG_{\mathrm{old}}\in L^r(\R^+;L^p(\Omega))$ and $\overline \bu\in L^r(0,T;D(A_p))$ then the problem~\eqref{eq:1338} admits a unique solution $\bG \in L^\infty(\R^+ \! \times \! (0,T);W^{1,p}(\Omega))$ such that $\partial_s \bG , \partial_t \bG \in L^\infty(\R^+;L^r(0,T;L^p(\Omega)))$.
Moreover, this solution satisfies
\begin{equation}\label{eq:1112a}
\begin{aligned}
\|\bG\|_{L^\infty(\R^+\! \times (0,T);W^{1,p}(\Omega))} + & \|\partial_s \bG\|_{L^\infty(\R^+;L^r(0,T;L^p(\Omega)))} + \|\partial_t \bG\|_{L^\infty(\R^+;L^r(0,T;L^p(\Omega)))} \\
& \qquad \leq C_2 \big( 1+ \|\nabla \overline \bu\|_{L^r(0,T;L^p(\Omega))} \big) \mathrm{exp} \big( C_3 \|\nabla \overline \bu\|_{L^1(0,T;W^{1,p}(\Omega))} \big),
\end{aligned}
\end{equation}
where the constants $C_2$ and $C_3$ depend on $\Omega$, $p$, $r$, $\We$ and the norms of $\bG_{\mathrm{old}}$ and $\partial_s\bG_{\mathrm{old}}$.
Their expressions will be precised in the proof.
\end{lemma}
The existence of a unique solution to~\eqref{eq:1338} follows from the application of the method of characteristics.
Some details are given in~\cite{Fernandez-Guillen-Ortega} (Appendix p.~26) about the equation~\eqref{eq:1650}, that is using the function~$\bF$.
In practice, the following estimates will be made on regular solution~$\bG_n$ which approaches the solution~$\bG$ when a regular velocity field~$\bu_n$ approaches the velocity~$\bu$.
The regularity of these solutions~$\bG_n$ with respect to~$t$ and~$s$ comes from the Cauchy-Lipschitz theorem.
For sake of simplicity, we omit the indexes~''$n$''.
In the following proof, we refer to~\cite{Fernandez-Guillen-Ortega} for the passage to the limit $n\to +\infty$.
The rest of the proof of the lemma~\ref{lem:1117} is split into three parts: in the first one (see the subsection~\ref{part:estimates-1}) we obtain a first estimate concerning the regularity of~$\bG$, and in the subsection~\ref{part:estimates-3} we obtain the estimate for~$\partial_t \bG$.
This estimate requires an estimate on~$\partial_s \bG$, which is given in the subsection~\ref{part:estimates-2}.
Note that we strongly use a Gronwall type lemma whose the proof is given in the Appendix~\ref{part:appendix1}.

\subsubsection{Estimate for the deformation gradient~$\bG$}\label{part:estimates-1}

Let $p>d$ and take the scalar product of the equation~\eqref{eq:1338} by $|\bG|^{p-2} \bG$.
We deduce
\begin{equation}
\frac{1}{p} \partial_t \big( |\bG|^p \big) + \frac{1}{\We\,p} \partial_s \big( |\bG|^p \big) + \frac{1}{p} \overline \bu \cdot \nabla \big( |\bG|^p \big) = |\bG|^{p-2} (\bG\cdot \nabla \overline \bu) : \bG.
\end{equation}
Integrating for $\bx\in \Omega$, due to the incompressible condition $\div\, \overline \bu = 0$ and the homogeneous boundary Dirichlet condition for the velocity, we obtain
\begin{equation}
\begin{aligned}
\partial_t \big( \|\bG\|_{L^p(\Omega)}^p \big) + \frac{1}{\We} \partial_s \big( \|\bG\|_{L^p(\Omega)}^p \big) 
= p \int_\Omega |\bG|^{p-2} (\bG\cdot \nabla \overline \bu) : \bG
\leq p \|\nabla \overline \bu\|_{L^\infty(\Omega)} \|\bG\|_{L^p(\Omega)}^p.
\end{aligned}
\end{equation}
We next use the continuous injection $W^{1,p}(\Omega) \hookrightarrow L^\infty(\Omega)$, holds for $p>d$ and making appear a constant $C_0=C_0(\Omega,p)$:
\begin{equation}\label{eq:1712}
\begin{aligned}
\partial_t \big( \|\bG\|_{L^p(\Omega)}^p \big) + \frac{1}{\We} \partial_s \big( \|\bG\|_{L^p(\Omega)}^p \big) 
\leq
p\, C_0 \|\nabla \overline \bu\|_{W^{1,p}(\Omega)} \|\bG\|_{L^p(\Omega)}^p.
\end{aligned}
\end{equation}
Now, we take the spatial gradient in~\eqref{eq:1338} and compute the scalar product of both sides of the resulting equation with $|\nabla \bG|^{p-2} \nabla \bG$ (we will note that this is a scalar product on the $3$-tensor, defined by $A::B = A_{i,j,k}B_{i,j,k}$).
After integrating for $\bx\in \Omega$ we obtain
\begin{equation}
\partial_t \big( \|\nabla \bG\|_{L^p(\Omega)}^p \big) + \frac{1}{\We} \partial_s \big( \|\nabla \bG\|_{L^p(\Omega)}^p \big)
\leq 2p \int_\Omega |\nabla \bG|^{p} |\nabla \overline \bu|
+ p \int_\Omega |\bG| |\nabla \bG|^{p-1} |\nabla^2\overline \bu|.
\end{equation}
Using the Hölder inequality, we have
\begin{equation}
\partial_t \big( \|\nabla \bG\|_{L^p(\Omega)}^p \big) + \frac{1}{\We} \partial_s \big( \|\nabla \bG\|_{L^p(\Omega)}^p \big)
\leq
2p \|\nabla \overline \bu\|_{L^\infty(\Omega)} \|\nabla \bG\|_{L^p(\Omega)}^p + p \|\bG\|_{L^\infty(\Omega)} \|\nabla \bG\|_{L^p(\Omega)}^{p-1}  \|\nabla^2 \overline \bu\|_{L^p(\Omega)} .
\end{equation}
For $p>d$, using the continuous injection $W^{1,p}(\Omega) \hookrightarrow L^\infty(\Omega)$ again, we deduce
\begin{equation}\label{eq:1713}
\partial_t \big( \|\nabla \bG\|_{L^p(\Omega)}^p \big) + \frac{1}{\We} \partial_s \big( \|\nabla \bG\|_{L^p(\Omega)}^p \big)
\leq 3p\, C_0 \|\nabla \overline \bu\|_{W^{1,p}(\Omega)} \|\bG\|_{W^{1,p}(\Omega)}^p.
\end{equation}
Adding this estimate~\eqref{eq:1713} with the estimate~\eqref{eq:1712}, we deduce
\begin{equation}
\partial_t \big( \|\bG\|_{W^{1,p}(\Omega)} \big) + \frac{1}{\We} \partial_s \big( \|\bG\|_{W^{1,p}(\Omega)} \big)
\leq 3C_0 \|\nabla \overline \bu\|_{W^{1,p}(\Omega)} \|\bG\|_{W^{1,p}(\Omega)}.
\end{equation}
Using the initial conditions we have $\|\bG\|_{W^{1,p}(\Omega)}\big|_{s=0} = \sqrt{d}\, |\Omega|^{\frac{1}{p}}$ and $\|\bG\|_{W^{1,p}(\Omega)}\big|_{t=0} = \|\bG_{\mathrm{old}}\|_{W^{1,p}(\Omega)}$, the Gronwall type lemma (see the Appendix~\ref{part:appendix1}) implies that for all $(s,t)\in \R^+\!\times \! (0,T)$ we have
\begin{equation}\label{eq:1604}
\|\bG(s,t,\cdot)\|_{W^{1,p}(\Omega)} \leq \zeta(s,t) \mathrm{exp} \Big( 3C_0 \int_0^t \|\nabla \overline \bu\|_{W^{1,p}(\Omega)}\Big),
\end{equation}
where $\dsp \zeta(s,t) = \left\{
\begin{aligned}
\|\bG_{\mathrm{old}}\|_{W^{1,p}(\Omega)}\big( s-\frac{t}{\We}\big) \quad & \text{if $t\leq \We\, s$}, \\
\sqrt{d}\, |\Omega|^{\frac{1}{p}} \hspace{2cm} & \text{if $t > \We\, s$}.
\end{aligned}
\right.
$
\par\noindent
The assumption $\bG_{\mathrm{old}}\in L^\infty(\R^+;W^{1,p}(\Omega))$ implies $\zeta \in L^\infty(\R^+ \! \times \! (0,T))$ with 
\begin{equation}
\|\zeta\|_{L^\infty(\R^+\times (0,T))}\leq \max\big\{\|\bG_{\mathrm{old}}\|_{L^\infty(\R^+;W^{1,p}(\Omega))} , \sqrt{d}\, |\Omega|^{\frac{1}{p}} \big\}.
\end{equation}
The relation~\eqref{eq:1604} now reads
\begin{equation}\label{eq:1610}
\begin{aligned}
& \|\bG\|_{L^\infty(\R^+\times (0,T);W^{1,p}(\Omega))}
\leq
\|\zeta\|_{L^\infty(\R^+\times (0,T))} \, \mathrm{exp}\big( 3C_0 \|\nabla \overline \bu\|_{L^1(0,T;W^{1,p}(\Omega))} \big).\\
\end{aligned}
\end{equation}

\subsubsection{Estimate for the age derivate~$\partial_s\bG$}\label{part:estimates-2}

We first remark that the derivative $\bG'=\partial_s \bG$ exactly satisfies the same PDE that $\bG$ (see the first equation of~\eqref{eq:1338}; that is due to the fact that $\overline \bu$ does not depend on the variable~$s$).
We then deduce the same kind of estimate that~\eqref{eq:1712}:
\begin{equation}\label{eq:2145}
\partial_t \big( \|\bG'\|_{L^p(\Omega)} \big) + \frac{1}{\We} \partial_s \big( \|\bG'\|_{L^p(\Omega)} \big)
\leq C_0 \|\nabla \overline \bu\|_{W^{1,p}(\Omega)} \|\bG'\|_{L^p(\Omega)}.
\end{equation}
But the initial conditions differ as follows: $\bG'|_{t=0}=\partial_s\bG_{\mathrm{old}}$ and $\bG'|_{s=0}=\We\, \nabla \overline \bu$.
This last condition is obtained using $s=0$ in the equation~\eqref{eq:1338}.
Note that this result is valid because we are working on regular solutions $\bG_n$ (see the introduction of this proof) such that $\partial_t \bG_n$ is continuous at $s=0$.
From the lemma~\ref{lem:gronwall} given in the Appendix~\ref{part:appendix1} we obtain for all $(s,t)\in \R^+\!\times \! (0,T)$ the estimate
\begin{equation}\label{eq:1605}
\|\bG'(s,t,\cdot)\|_{L^p(\Omega)}
\leq
\zeta'(s,t) \, \mathrm{exp} \Big( C_0 \int_0^t \|\nabla \overline \bu\|_{W^{1,p}(\Omega)}\Big),
\end{equation}
where $ \dsp \zeta'(s,t)=\left\{
\begin{aligned}
\|\partial_s\bG_{\mathrm{old}}\|_{L^p(\Omega)}\big( s-\frac{t}{\We} \big) \quad & \text{if $t\leq \We\, s$},\\
\We \, \|\nabla \overline \bu\|_{L^p(\Omega)}\big( t-\We\, s \big) \quad & \text{if $t > \We\, s$}.\\
\end{aligned}
\right.
$
\par\noindent
For each $s\geq 0$, we estimate the $L^r(0,T)$-norm of the function $t\mapsto \zeta'(s,t)$ as follows: if $T\leq \We\, s$ then
\begin{equation}
\begin{aligned}
\int_0^T \zeta'(s,t)^r \, \mathrm dt 
& = \int_0^T \|\partial_s\bG_{\mathrm{old}}\|_{L^p(\Omega)}^r\big( s-\frac{t}{\We} \big) \, \mathrm dt \\ 
& = \We \int_{s-\frac{T}{\We}}^s \|\partial_s\bG_{\mathrm{old}}\|_{L^p(\Omega)}^r(t') \, \mathrm dt' \\
& \leq \We \, \|\partial_s\bG_{\mathrm{old}}\|_{L^r(\R^+;L^p(\Omega))}^r.
\end{aligned}
\end{equation}
If $T > \We\, s$ then
\begin{equation}
\begin{aligned}
\int_0^T \zeta'(s,t)^r \, \mathrm dt 
& = \int_0^{\We\, s} \|\partial_s\bG_{\mathrm{old}}\|_{L^p(\Omega)}^r\big( s-\frac{t}{\We} \big) \, \mathrm dt + \We^r \int_{\We\, s}^T \|\nabla \overline \bu\|_{L^p(\Omega)}^r\big( t-\We\, s \big) \, \mathrm dt \\ 
& = \We \int_0^s \|\partial_s\bG_{\mathrm{old}}\|_{L^p(\Omega)}^r(t') \, \mathrm dt' + \We^r \int_0^{T-\We\, s} \|\nabla \overline \bu\|_{L^p(\Omega)}^r(t') \, \mathrm dt' \\
& \leq \We \, \|\partial_s\bG_{\mathrm{old}}\|_{L^r(\R^+;L^p(\Omega))}^r + \We^r \, \|\nabla \overline \bu\|_{L^r(0,T;L^p(\Omega))}^r.
\end{aligned}
\end{equation}
Finally, we obtain $\zeta'\in L^\infty(\R^+;L^r(0,T))$ with
\begin{equation}
\|\zeta'\|_{L^\infty(\R^+;L^r(0,T))} \leq \We^{\frac{1}{r}} \|\partial_s \bG_{\mathrm{old}}\|_{L^r(\R^+;L^p(\Omega))} + \We \|\nabla \overline \bu\|_{L^r(0,T;L^p(\Omega))}.
\end{equation}
The relation~\eqref{eq:1605} now reads
\begin{equation}\label{eq:1611}
\|\bG'\|_{L^\infty(\R^+;L^r(0,T;L^p(\Omega)))} \leq \|\zeta'\|_{L^\infty(\R^+;L^r(0,T))} \, \mathrm{exp} \big( C_0 \|\nabla \overline \bu\|_{L^1(0,T;W^{1,p}(\Omega))} \big).
\end{equation}

\subsubsection{Estimate for the time derivate~$\partial_t\bG$}\label{part:estimates-3}

Isolating the term $\partial_t\bG$ in the equation~\eqref{eq:1338} we have
\begin{equation}
\begin{aligned}
\|\partial_t \bG\|_{L^p(\Omega)}
& \leq
\frac{1}{\We} \|\bG'\|_{L^p(\Omega)} + \|\overline \bu\|_{L^\infty(\Omega)} \|\nabla \bG\|_{L^p(\Omega)} + \|\bG\|_{L^\infty(\Omega)} \|\nabla \overline \bu\|_{L^p(\Omega)} \\
& \leq
\frac{1}{\We} \|\bG'\|_{L^p(\Omega)} + C_0 \|\overline \bu\|_{W^{1,p}(\Omega)} \|\nabla \bG\|_{L^p(\Omega)} + C_0 \|\bG\|_{W^{1,p}(\Omega)} \|\nabla \overline \bu\|_{L^p(\Omega)}.
\end{aligned}
\end{equation}
Introducing the Poincar\'e inequality with a constant $C_P=C_P(\Omega,p)$, holds since $\overline \bu$ vanishes on the boundary of the domain, we obtain
\begin{equation}
\|\partial_t \bG\|_{L^p(\Omega)}
\leq
\frac{1}{\We} \|\bG'\|_{L^p(\Omega)} + C_0 (1+C_P) \|\nabla \overline \bu\|_{L^p(\Omega)} \|\bG\|_{W^{1,p}(\Omega)}.
\end{equation}
Taking the $L^r(0,T)$-norm for the variable~$t$, and next the $L^{\infty}(\R^+)$-norm for the variable~$s$, we obtain
\begin{equation}
\begin{aligned}
\|\partial_t \bG\|_{L^\infty(\R^+;L^r(0,T;L^p(\Omega)))}
\leq
& \frac{1}{\We} \|\bG'\|_{L^\infty(\R^+;L^r(0,T;L^p(\Omega)))} \\
& \quad + C_0 (1+C_P) \|\nabla \overline \bu\|_{L^r(0,T;L^p(\Omega))} \|\bG\|_{L^{\infty}(\R^+\times (0,T);W^{1,p}(\Omega))}.
\end{aligned}
\end{equation}
Using the previous estimates~\eqref{eq:1610} and~\eqref{eq:1611}, we deduce the result announced in the lemma~\ref{lem:1117}.

\subsection{Estimates for the stress tensor~$\btau$}

\begin{lemma}\label{lem:1118}
Let $T>0$, $r\in ]1,+\infty[$ and $p\in]d,+\infty[$.
\par\noindent
If $\overline \bG \in L^\infty(\R^+ \! \times \! (0,T);W^{1,p}(\Omega))$ and $\partial_t \overline \bG \in L^\infty(\R^+;L^r(0,T;L^p(\Omega)))$ then the stress tensor~$\btau$ defined by the integral relation~\eqref{eq:1339} belongs to $L^\infty(0,T;W^{1,p}(\Omega))$ and its time derivative $\partial_t \btau$ belongs to $L^s(0,T;L^p(\Omega))$.
Moreover, there exists a continuous increasing function $F_0:\R^+\mapsto \R^+$ such that
\begin{equation}\label{eq:1113}
\|\btau\|_{L^\infty(0,T;W^{1,p}(\Omega))} + \|\partial_t \btau\|_{L^r(0,T;L^p(\Omega))} \leq \frac{\omega}{\We} F_0 \big( \|\overline \bG\|_{L^\infty(\R^+\times (0,T);W^{1,p}(\Omega))} + \|\partial_t \overline \bG\|_{L^\infty(\R^+;L^r(0,T;L^p(\Omega)))} \big).
\end{equation}
The function $F_0$ depends on $\Omega$, $p$ and on the function growth of the function~$\mathcal S$.
\end{lemma}

\proof
Since the function $\mathcal S$ is of class $\mathcal C^1$, we can introduce the following continuous and non-decreasing real functions
\begin{equation}\label{eq:1626}
\begin{aligned}
& \mathcal S_0:c\in \R^+ \longmapsto \max_{|\bG|\leq c}|\mathcal S(\bG)|\in \R^+, \\
& \mathcal S_1:c\in \R^+ \longmapsto \max_{|\bG|\leq c}|\mathcal S'(\bG)|\in \R^+,
\end{aligned}
\end{equation}
where the derivative $\mathcal S'(\bG)$ denotes the $4$-tensor whose the coefficient $(i,j,k,\ell)$, denoted $\partial_{(ij)}\mathcal S(\overline \bG)_{k\ell}$, is the derivative of $\big( \mathcal S(\bG)\big)_{k\ell}$ with respect to the tensor $\boldsymbol{E}_{ij}$ of the canonical basis of the space~$\mathcal L(\R^d)$ of real $d\times d$ matrices, see the Appendix~\ref{part:appendix0}.
\par\noindent
Due to the continuous injection $W^{1,p}(\Omega) \hookrightarrow L^\infty(\Omega)$, the function $\overline \bG$ introduced in the hypothesis of the lemma is bounded in $\R^+ \! \times \! [0,T] \times \Omega$ and we have
\begin{equation}
\begin{aligned}
\|\mathcal S(\overline \bG)\|_{L^\infty(\R^+\times (0,T) \times \Omega)} 
& \leq \mathcal S_0(\|\overline \bG\|_{L^\infty(\R^+\times (0,T) \times \Omega)}) \\
& \leq \mathcal S_0(C_0\, \|\overline \bG\|_{L^\infty(\R^+\times (0,T);W^{1,p}(\Omega))}).
\end{aligned}
\end{equation}
To simplify, we note $\overline{c}:=C_0\, \|\overline \bG\|_{L^\infty(\R^+\times (0,T);W^{1,p}(\Omega))}$.
In the same way the function $\mathcal S'(\overline \bG)$ is bounded by the real~$\mathcal S_1(\overline{c})$.\\[0.3cm]
{\bf $\boldsymbol{L^p}$-norm for $\btau$ --}
We easily have the following bound for the stress tensor $\btau$ given by the formula~\eqref{eq:1339}: $| \btau(t,\bx) | \leq \frac{\omega}{\We} \mathcal S_0(\overline{c})$ for a.~e. $(t,\bx)\in (0,T)\!\times \!\Omega$.
We note that we used $\int_0^\infty m = 1$.
We deduce in particular that
\begin{equation}\label{eq:1601}
\| \btau \|_{L^\infty(0,T;L^p(\Omega))} \leq |\Omega|^{\frac{1}{p}} \frac{\omega}{\We} \mathcal S_0(\overline{c}).
\end{equation}
{\bf $\boldsymbol{W^{1,p}}$-norm for $\btau$ --}
Taking the spatial gradient of the expression~\eqref{eq:1339} given the stress tensor we obtain
\begin{equation}
\nabla \btau(t,\bx) = 
\frac{\omega}{\We} \int_0^{+\infty} m(s) \nabla \overline \bG(s,t,\bx) : \mathcal S'(\overline \bG(s,t,\bx))\, \mathrm ds.
\end{equation}
The meaning of the symbols here is the following.
Component by component, the equality above written
\begin{equation}
\big[ \nabla \btau(t,\bx) \big]_{ijk} = \partial_i \btau_{jk}(t,\bx) = 
\frac{\omega}{\We} \int_0^{+\infty} m(s) \sum_{\ell,m} \partial_i \overline \bG_{\ell m}(s,t,\bx)\, \partial_{(\ell m)}\mathcal S\big(\overline \bG(s,t,\bx)\big)_{jk}\, \mathrm ds.
\end{equation}
Using the Hölder inequality, the $L^\infty$-bound on $\mathcal S'(\overline \bG)$ and the fact that $\dsp \int_0^\infty m = 1$, we obtain
\begin{equation}
  \begin{aligned}
| \nabla \btau(t,\bx) |^p 
& = \frac{\omega^p}{\We^p} \Big| \int_0^\infty  m(s)^{\frac{1}{p}} \nabla \overline \bG(s,t,\bx) : m(s)^{1-\frac{1}{p}} \mathcal S'\big(\overline \bG(s,t,\bx)\big) \, \mathrm ds \Big|^p \\
& \leq \frac{\omega^p}{\We^p} \mathcal S_1(\overline{c})^p \int_0^\infty m(s) |\nabla \overline \bG(s,t,\bx)|^p \, \mathrm ds.
  \end{aligned}
\end{equation}
Integrating for $\bx\in \Omega$, we obtain
\begin{equation}
\| \nabla \btau(t,\cdot) \|_{L^p(\Omega)}^p 
\leq
\frac{\omega^p}{\We^p} \mathcal S_1(\overline{c})^p  \int_0^\infty m(s) \|\nabla \overline \bG(s,t,\cdot)\|_{L^p(\Omega)}^p \, \mathrm ds.
\end{equation}
Due to the definition of the bound~$\overline c$, this implies
\begin{equation}\label{eq:1602}
\| \nabla \btau \|_{L^\infty(0,T;L^p(\Omega))}
\leq
\frac{\omega}{\We} \mathcal S_1(\overline c) \, \frac{\overline c}{C_0}.
\end{equation}
{\bf $\boldsymbol{L^p}$-norm for $\boldsymbol{\partial_t} \btau$ --}
Similarly, we obtain a bound for $\partial_t \btau$ in $L^p(\Omega)$: we have
\begin{equation}\label{eq:2221}
\partial_t \btau(t,\bx) = 
\frac{\omega}{\We} \int_0^{+\infty} m(s) \partial_t \overline \bG(s,t,\bx) : \mathcal S'\big(\overline \bG(s,t,\bx)\big)\, \mathrm ds.
\end{equation}
Using the Hölder inequality to control the quantity $|\partial_t \btau(t,\bx)|^p$, and next an integration for $\bx\in \Omega$, we obtain
\begin{equation}
\| \partial_t \btau(t,\cdot) \|_{L^p(\Omega)}^p 
\leq
\frac{\omega^p}{\We^p} \mathcal S_1(\overline c)^p \int_0^\infty m(s) \|\partial_t \overline \bG(s,t,\cdot)\|_{L^p(\Omega)}^p \, \mathrm ds.
\end{equation}
Due to the assumption on $\partial_t \bG$ which implies that $\widetilde c:= \|\partial_t \overline \bG\|_{L^\infty(\R^+;L^r(0,T;L^p(\Omega)))} < +\infty$, we deduce
\begin{equation}\label{eq:1603}
\| \partial_t \btau \|_{L^r(0,T;L^p(\Omega))}
\leq
\frac{\omega}{\We}  \mathcal S_1(\overline c) \, \widetilde c.
\end{equation}
The estimates~\eqref{eq:1601}, \eqref{eq:1602} and~\eqref{eq:1603} show that $\btau$ and $\partial_t \btau$ are bounded respectively in $L^\infty(0,T;W^{1,p}(\Omega))$ and in $L^r(0,T;L^p(\Omega))$, and that these bounds continuously depend on~$\overline c$ and~$\widetilde c$, and increase in both variables.
\cqfd

\subsection{Proof of the theorem~\ref{th:local}}\label{part:estimates-4}

For any $T>0$ we introduce the Banach space 
$$
\mathscr{B}(T) = L^r(0,T;W^{1,p}_0(\Omega)) \times \mathcal C(\R^+ \! \times \! [0,T];L^p(\Omega)) \times \mathcal C([0,T];L^p(\Omega))
$$
and for any $R_1>0$, $R_2>0$ and $R_3>0$ the subset
\begin{equation}
\begin{aligned}
& \hspace{-0.3cm} \mathscr{H}(T,R_1,R_2,R_3)
=
\Big\{ (\overline{\bu} , \overline \bG, \overline{\btau}) \in \mathscr{B}(T) ~ ; ~
\overline \bu \in L^r(0,T;D(A_p)),
\quad \partial_t \overline \bu \in L^r(0,T;H_p), \\
& \overline \bG \in L^\infty(\R^+\!\times\! (0,T);W^{1,p}(\Omega)),
\quad  \partial_s \overline \bG, ~\partial_t \overline \bG \in L^\infty(\R^+;L^r(0,T;L^p(\Omega))),\\
& \overline \btau \in L^\infty(0,T;W^{1,p}(\Omega)),
\quad  \partial_t \overline \btau \in L^r(0,T;L^p(\Omega)),\\
& \overline \bu|_{t=0} = \bu_0,
\quad \overline \bG|_{t=0} = \bG_{\mathrm{old}},
\quad \overline \bG|_{s=0} = \bdelta, \\
& \|\overline \bu\|_{L^r(0,T;W^{2,p}(\Omega))} + \| \partial_t \overline \bu\|_{L^r(0,T;L^p(\Omega))} \leq R_1, \\
& \|\overline \bG\|_{L^\infty(\R^+\times (0,T);W^{1,p}(\Omega))} + \|\partial_s \overline \bG\|_{L^\infty(\R^+;L^r(0,T;L^p(\Omega)))} + \|\partial_t \overline \bG\|_{L^\infty(\R^+;L^r(0,T;L^p(\Omega)))} \leq R_2, \\
& \|\overline \btau\|_{L^\infty(0,T;W^{1,p}(\Omega))} + \|\partial_t \overline \btau\|_{L^r(0,T;L^p(\Omega))} \leq R_3 ~ \Big\}.
\end{aligned}
\end{equation}
\begin{remark}
Such a set is non empty, for instance if $R_1$ and~$R_2$ are large enough.
More precisely, if
\begin{equation}\label{eq:2128}
R_1 \geq \frac{C_1}{1-\omega}\|\bu_0\|_{D^r_p(\Omega)}
\qquad \text{and} \qquad
R_2 \geq \|\bG_{\mathrm{old}}\|_{L^\infty(\R^+;W^{1,p}(\Omega))}
\end{equation}
then for any $T>0$ and any $R_3>0$ we can build a velocity field $\bu^\star$ such that $(\bu^\star,\bG_{\mathrm{old}},\b0) \in \mathscr{H}(T,R_1,R_2,R_3)$, see an example of construction in~\cite{Fernandez-Guillen-Ortega, Guillope-Saut3}.
\end{remark}
\begin{remark}
If $(\bu,\bG,\btau)\in \mathscr{H}(T,R_1,R_2,R_3)$ for some $T$, $R_1$, $R_2$ and~$R_3$ then the velocity field~$\bu$ and the tensor~$\bG$ are continuous with respect to the time~$t$ and the age~$s$.
In fact, these continuity properties follow from the Sobolev injections of kind $W^{1,\alpha}(0,A;X) \subset \mathcal C([0,A];X)$, holds for $\alpha>1$.
Moreover, they make sense of the initial conditions $\bu|_{t=0} = \bu_0$, $\bG|_{t=0} = \bG_{\mathrm{old}}$ and~$\bG|_{s=0} = \bdelta$.
\end{remark}
We consider the mapping
\begin{equation}
\begin{aligned}
\Phi ~:~ \mathscr{H}(T,R_1,R_2,R_3) & \longrightarrow ~ \mathscr{B}(T) \\
(\overline{\bu} , \overline{\bG} , \overline{\btau}) ~~~~~ & \longmapsto ~(\bu,\bG,\btau),
\end{aligned}
\end{equation}
where $\bu$ is the unique solution of the Stokes problem~\eqref{eq:1337} with
\begin{equation}\label{eq:g}
\overline \bg = -\Re\, \overline{\bu} \cdot \nabla \overline{\bu} + \div\, \overline{\btau} + \boldsymbol{f};
\end{equation}
where $\bG$ solves the problem~\eqref{eq:1338} and where $\btau$ is given by the integral formula~\eqref{eq:1339}.
The goal of this proof is to show that the application~$\Phi$ has a fixed point.
For this we first prove that $\Phi$ leaves a set $\mathscr{H}(T,R_1,R_2,R_3)$ invariant (for a ``good'' choice of $T$, $R_1$, $R_2$ and $R_3$).
\par\noindent
Let $T>0$, $R_1>0$, $R_2>0$, $R_3>0$  and $(\overline{\bu} , \overline{\bG} , \overline{\btau}) \in \mathscr{H}(T,R_1,R_2,R_3)$.
If we denote by $(\bu,\bG,\btau) = \Phi(\overline{\bu} , \overline{\bG} , \overline{\btau})$, we will show that the previous lemmas imply estimates of $(\bu,\bG,\btau)$ with respect to the norms of $(\overline{\bu} , \overline{\bG} , \overline{\btau})$, that is with respect to $(T,R_1,R_2,R_3)$.\\[0.3cm]
{\bf Velocity estimate --}
From the lemma~\ref{lem:stokes} we can estimate~$\bu$ and~$\partial_t \bu$ using the norm $\|\bg\|_{L^r(0,T;L^p(\Omega))}$.
For the source term~$\bg$ given by the relation~\eqref{eq:g}, we have
\begin{equation}
\|\bg\|_{L^r(0,T;L^p(\Omega))} \leq \Re \underbrace{\|\overline \bu \cdot \nabla \overline \bu\|_{L^r(0,T;L^p(\Omega))}}_{T_1} + \underbrace{\|\overline \btau\|_{L^r(0,T;W^{1,p}(\Omega))}}_{T_2} + \|\boldsymbol{f}\|_{L^r(0,T;L^p(\Omega))}.
\end{equation}
Since we have the bound~$R_3$ on~$\overline{\btau}$ in $L^\infty(0,T;W^{1,p}(\Omega))$, the term $T_2$ satisfies $T_2 \leq T^{\frac{1}{r}} R_3$.
The bilinear term $T_1$ is more difficult to estimate.
We follow the ideas of~\cite{Fernandez-Guillen-Ortega} and we generalize their result to the $d$-dimensional case (the paper~\cite{Fernandez-Guillen-Ortega} only deals with the case $d=3$):
\begin{equation}\label{eq:1329}
\begin{aligned}
T_1
& \leq \|\overline \bu\|_{L^{2r}(0,T;L^\infty(\Omega))} \|\nabla \overline \bu\|_{L^{2r}(0,T;L^p(\Omega))} \\
& \leq T^{\frac{p-d}{2rp}} \|\overline \bu\|_{L^{\frac{2rp}{d}}(0,T;L^\infty(\Omega))} \|\overline \bu\|_{L^{2r}(0,T;W^{1,p}(\Omega))}.
\end{aligned}
\end{equation}
But we have the following estimate (see~\cite{Friedman}):
\begin{equation}
\|\overline \bu\|_{L^\infty(\Omega)} \leq C \|\overline \bu\|_{L^p(\Omega)}^{\frac{p-d}{p}} \|\overline \bu\|_{W^{1,p}(\Omega)}^{\frac{d}{p}},
\end{equation}
which, after integrating with respect to time, implies
\begin{equation}\label{eq:1340}
\|\overline \bu\|_{L^{\frac{2pr}{d}}(0,T;L^\infty(\Omega))} \leq C \|\overline \bu\|_{L^\infty(0,T;L^p(\Omega))}^{\frac{p-d}{p}} \|\overline \bu\|_{L^{2r}(0,T;W^{1,p}(\Omega))}^{\frac{d}{p}}.
\end{equation}
Note that the constant~$C$ introduced here only depends on~$\Omega$, $p$ and $d$.
Moreover, by interpolation, we have
\begin{equation}\label{eq:1341}
\|\overline \bu\|_{L^{2r}(0,T;W^{1,p}(\Omega))} \leq \|\overline \bu\|_{L^\infty(0,T;L^p(\Omega))}^{\frac{1}{2}} \|\overline \bu\|_{L^r(0,T;W^{2,p}(\Omega))}^{\frac{1}{2}}.
\end{equation}
Using~\eqref{eq:1340} and~\eqref{eq:1341}, the estimate~\eqref{eq:1329} now reads
\begin{equation}\label{eq:1330}
\begin{aligned}
T_1
& \leq C\, T^{\frac{p-d}{2rp}} \|\overline \bu\|_{L^\infty(0,T;L^p(\Omega))}^{\frac{3p-d}{2p}} \|\overline \bu\|_{L^r(0,T;W^{2,p}(\Omega))}^{\frac{p+d}{2p}}.
\end{aligned}
\end{equation}
Finally, we use $\overline \bu(t,\bx) = \bu_0(\bx) + \int_0^t \partial_t \overline \bu$ to obtain
\begin{equation}
\begin{aligned}
\|\overline \bu\|_{L^\infty(0,T;L^p(\Omega))}
& \leq \|\bu_0\|_{L^p(\Omega)} + \|\partial_t \overline \bu\|_{L^1(0,T;L^p(\Omega))} \\
& \leq \|\bu_0\|_{L^p(\Omega)} + T^{1-\frac{1}{r}} \|\partial_t \overline \bu\|_{L^r(0,T;L^p(\Omega))}.
\end{aligned}
\end{equation}
Using the bound~$R_1$ for $\overline \bu$ and $\partial_t \overline \bu$ given in the definition of the subspace~$\mathscr{H}(T,R_1,R_2,R_3)$, the estimate~\eqref{eq:1330} becomes
\begin{equation}\label{eq:1331}
T_1 \leq C\, T^{\frac{p-d}{2rp}} R_1^{\frac{p+d}{2p}} \|\bu_0\|_{L^p(\Omega)}^{\frac{3p-d}{2p}}
+ C\, T^{\frac{3p-d}{2p}-\frac{1}{r}} R_1^2.
\end{equation}
We now use this bound to control the source term~$\bg$.
The lemma~\ref{lem:stokes} implies:
\begin{equation}\label{eq:1499}
\begin{aligned}
\|\bu\|_{L^r(0,T;W^{2,p}(\Omega))} + \| \partial_t \bu\|_{L^r(0,T;L^p(\Omega))}
\leq
\frac{C_1}{1-\omega} \Big( & \|\bu_0\|_{D^r_p(\Omega)} + \|\boldsymbol{f}\|_{L^r(0,T;L^p(\Omega))} \\
& + \Re \, C\, T^\alpha R_1^\beta \|\bu_0\|_{D^r_p(\Omega)}^\gamma + \Re \, C\, T^\delta R_1^2 + T^{\frac{1}{r}} R_3  \Big),
\end{aligned}
\end{equation}
where the assumptions $p>d$ and $r>1$ imply that $\alpha$, $\beta$, $\gamma$ and $\delta$ are positive numbers.
This estimate~\eqref{eq:1499} can be rewrite as
\begin{equation}\label{eq:1500}
\|\bu\|_{L^r(0,T;W^{2,p}(\Omega))} + \| \partial_t \bu\|_{L^r(0,T;L^p(\Omega))}
\leq
\widetilde{C_1} \Big( 1 + T^{\frac{1}{r}}\, R_3 + K(T,R_1)\Big),
\end{equation}
where $\widetilde{C_1}$ may also depends on $\omega$ and on the norm of $\bu_0$ and~$\boldsymbol{f}$ in their spaces.
It is important to notice that for each $R_1>0$ we have $\dsp \lim_{T\to 0} K(T,R_1) = 0$, and that $\widetilde{C_1}\geq \frac{C_1}{1-\omega} \|\bu_0\|_{D^r_p(\Omega)}$.\\[0.3cm]
{\bf Deformation gradient estimate --}
From the lemma~\ref{lem:1117}, we have
\begin{equation}\label{eq:1112}
\begin{aligned}
\|\bG\|_{L^\infty(\R^+\times (0,T);W^{1,p}(\Omega))} + \|\partial_t \bG\|_{L^\infty(\R^+;L^r(0,T;L^p(\Omega)))}
\leq
C_2 \big( 1+ R_1 \big) \mathrm{exp} \big( C_3 T^{1-\frac{1}{r}} R_1 \big).
\end{aligned}
\end{equation}
{\bf Stress tensor estimate --}
From the lemma~\ref{lem:1118}, we have
\begin{equation}\label{eq:1114}
\begin{aligned}
\|\btau\|_{L^\infty(0,T;W^{1,p}(\Omega))} + \|\partial_t \btau\|_{L^r(0,T;L^p(\Omega))}
\leq
\frac{\omega}{\We} F_0(R_2).
\end{aligned}
\end{equation}
{\bf $\Phi$-Invariant subset --}
If we then successively choose
\begin{equation}
\begin{aligned}
& R_1^\star = 2\widetilde{C_1}, \\
& R_2^\star = C_2 \big( 1+ R_1^\star \big) \mathrm{exp} \big( C_3 R_1^\star \big) + \|\bG_{\mathrm{old}}\|_{L^\infty(\R^+;W^{1,p}(\Omega))}, \\
& R_3^\star = \frac{\omega}{\We} F_0(R_2^\star),\\
& \text{and} \quad T_\star\leq 1 ~\text{small enough to have}~~T_\star^{\frac{1}{r}} R_3^\star + K(T_\star,R_1^\star) \leq 1,
\end{aligned}
\end{equation}
then we verify that $\mathscr{H}(T_\star,R_1^\star,R_2^\star,R_3^\star) \neq \emptyset$ (that is the inequalities~\eqref{eq:2128} hold).
For such a choice, the estimates~\eqref{eq:1500}, \eqref{eq:1112} and~\eqref{eq:1114} imply that $\Phi(\mathscr{H}(T_\star,R_1^\star,R_2^\star,R_3^\star)) \subset \mathscr{H}(T_\star,R_1^\star,R_2^\star,R_3^\star)$.
Moreover the function $\Phi$ is continuous and $\mathscr{H}(T_\star,R_1^\star,R_2^\star,R_3^\star)$ is a convex compact subset of $\mathscr{B}(T_\star)$, see~\cite{Guillope-Saut3} for similar properties.
We conclude the proof using the Schauder's theorem.
\cqfd

\section{Proof for the uniqueness result}\label{part:unic}

This section is devoted to the proof of the local existence theorem~\ref{th:unic}.
\par\noindent
As usual, we take the difference of the two solutions indexed by~$1$ and~$2$.
The vector $\bu=\bu_1-\bu_2$, the scalar $p=p_1-p_2$ and the tensors $\btau = \btau_1-\btau_2$, $\bG=\bG_1-\bG_2$ satisfy the following:
\begin{equation}\label{eq:1332}
\left\{
\begin{aligned}
& \Re ( \partial_t \bu + \bu_1\cdot \nabla \bu + \bu\cdot \nabla \bu_2 ) + \nabla p - (1-\omega) \Delta \bu = \div \, \btau, \\
& \div \, \bu = 0, \\
& \btau = \frac{\omega}{\We} \int_0^{+\infty} m(s) \, \Big[\mathcal S\big(\bG_1(s,\cdot,\cdot)\big) - \mathcal S\big(\bG_2(s,\cdot,\cdot)\big) \Big]\, \mathrm ds,\\
& \partial_t \bG + \frac{1}{\We} \partial_s \bG + \bu_1 \cdot \nabla \bG + \bu \cdot \nabla \bG_2 = \bG_1 \cdot \nabla \bu + \bG \cdot \nabla \bu_2,
\end{aligned}
\right.
\end{equation}
together with zero initial conditions $\bu\big|_{t=0} = \b0$ and $\bG\big|_{t=0} = \bG\big|_{s=0} = \b0$.
Note that the regularity of $\bG_i$ and the definition of the stress tensor $\btau_i$ implies that $\btau_i \in L^\infty(0,T;W^{1,d}(\Omega))$ (the proof is similar that those presented in the proof of the existence theorem~\ref{th:local}).
The uniqueness proof consists in demonstrate that $\bu = \b0$ and that $\bG=\btau=\b0$.
We will initially provide estimates on these three quantities.
\\[0.3cm]
{\bf Velocity estimate --}
Taking the scalar product of the first equation of the system~\eqref{eq:1332} by $\bu$ in $L^2(\Omega)$, we obtain
\begin{equation}
\frac{\Re}{2} d_t \big( \|\bu\|_{L^2(\Omega)}^2 \big) + (1-\omega) \|\nabla \bu\|_{L^2(\Omega)}^2
=
- \int_\Omega \btau \cdot \nabla \bu - \Re \int_\Omega (\bu \cdot \nabla \bu_2) \cdot \bu.
\end{equation}
From the Hölder inequality and the Young inequality, we obtain
\begin{equation}
\Re \, d_t \big( \|\bu\|_{L^2(\Omega)}^2 \big) + (1-\omega) \|\nabla \bu\|_{L^2(\Omega)}^2
\leq \frac{4}{1-\omega}\|\btau\|_{L^2(\Omega)}^2 + \Re \|\nabla \bu_2\|_{L^\infty(\Omega)} \|\bu\|_{L^2(\Omega)}^2.
\end{equation}
Introducing $Z(t)=\Re\,\|\bu\|_{L^2(\Omega)}^2$ and $C_Z(t) = \|\nabla \bu_2\|_{L^\infty(\Omega)} \in L^1(0,T)$, this estimate reads
\begin{equation}\label{eq:unicity1}
Z'(t) + (1-\omega) \|\nabla \bu\|_{L^2(\Omega)}^2
\leq
\frac{4}{1-\omega} \|\btau\|_{L^2(\Omega)}^2  + C_Z(t) Z(t).
\end{equation}
{\bf Stress tensor estimate --}
From the definition of the stress tensor~$\btau$ in the System~\eqref{eq:1332} and using the Cauchy-Schwarz inequality, we have, for all $(t,\bx)\in (0,T)\times \Omega$:
\begin{equation}
|\btau(t,\bx)|^2 \leq \frac{\omega^2}{\We^2} \int_0^{+\infty} m(s) \, \Big|\mathcal S\big(\bG_1(s,t,\bx)\big) - \mathcal S\big(\bG_2(s,t,\bx)\big) \Big|^2\, \mathrm ds.
\end{equation}
By assumption, the function $\mathcal S$ is of class~$\mathcal C^1$, so that~$\mathcal S'$ is bounded on each compact.
Since $\bG_i$, $i\in \{1,2\}$, belongs to $L^\infty(\R^+\times (0,T)\times \Omega)$ we deduce that there exists a constant $C'$, only depending on the norm $\|\bG_i\|_{L^\infty(\R^+\times (0,T)\times \Omega)}$ such that $|\mathcal S(\bG_1)-\mathcal S(\bG_2)|^2\leq C' |\bG_1-\bG_2|^2$ a.e. in $\R^+\! \times \! (0,T)\times \Omega$.
We deduce
\begin{equation}
|\btau(t,\bx)|^2 \leq \frac{C' \, \omega^2}{\We^2} \int_0^{+\infty} m(s) \, |\bG(s,t,\bx) |^2\, \mathrm ds.
\end{equation}
Integrating with respect to $\bx\in \Omega$ we obtain
\begin{equation}\label{eq:unicity2}
\|\btau\|_{L^2(\Omega)}^2 \leq \frac{C' \, \omega^2}{\We^2} Y(t),
\end{equation}
where we introduced $\dsp Y(t) = \int_0^{+\infty} \!\! m(s) \, \|\bG(s,t,\cdot)\|_{L^2(\Omega)}^2\, \mathrm ds$.\\[0.3cm]
{\bf Deformation gradient estimate --}
Taking the scalar product of the last equation of the system~\eqref{eq:1332} by $\bG$ in $L^2(\Omega)$, we obtain
\begin{equation}
\frac{1}{2} \partial_t \big( \|\bG\|_{L^2(\Omega)}^2 \big) + \frac{1}{2\We} \partial_s \big( \|\bG\|_{L^2(\Omega)}^2 \big) = \int_\Omega (\bG_1 \cdot \nabla \bu)\cdot \bG + \int_\Omega (\bG \cdot \nabla \bu_2)\cdot \bG - \int_\Omega (\bu \cdot \nabla \bG_2)\cdot \bG.
\end{equation}
Using the Hölder inequality, we have the estimate
\begin{equation}\label{eq:2205}
\begin{aligned}
\frac{1}{2} \partial_t \big( \|\bG\|_{L^2(\Omega)}^2 \big) + \frac{1}{2\We} \partial_s \big( \|\bG\|_{L^2(\Omega)}^2 \big)
\leq
& \|\bG_1\|_{L^\infty(\Omega)} \|\nabla \bu\|_{L^2(\Omega)} \|\bG\|_{L^2(\Omega)} \\
& + \|\nabla \bu_2\|_{L^\infty(\Omega)} \|\bG\|_{L^2(\Omega)}^2 \\
& + \|\bu\|_{L^{\frac{2d}{d-2}}(\Omega)} \|\nabla \bG_2\|_{L^d(\Omega)} \|\bG\|_{L^2(\Omega)}.
\end{aligned}
\end{equation}
Due to the Sobolev continuous injection $H^1(\Omega) \hookrightarrow L^{\frac{2d}{d-2}}(\Omega)$, the Poincar\'e inequality and the Young inequality, we obtain for all $\eps>0$:
\begin{equation}\label{eq:2320}
\partial_t \big( \|\bG\|_{L^2(\Omega)}^2 \big) + \frac{1}{\We} \partial_s \big( \|\bG\|_{L^2(\Omega)}^2 \big)
\leq
\eps \|\nabla \bu\|_{L^2(\Omega)}^2 + C_Y(t) \|\bG\|_{L^2(\Omega)}^2,
\end{equation}
where the function $C_Y(t) = \sup_{s\in \R^+}\Big\{ \frac{2}{\eps} \|\bG_1\|_{L^\infty(\Omega)}^2 + \|\nabla \bu_2\|_{L^\infty(\Omega)} + \frac{2\, C^2}{\eps} \, \|\nabla \bG_2\|_{L^d(\Omega)}^2\Big\}\in L^1(0,T)$, and where the constant~$C$ depends on~$\Omega$, $p$ and~$d$.
Multiplying this estimate~\eqref{eq:2320} by $m(s)$ and integrating for $s\in (0,+\infty)$ we obtain
\begin{equation}\label{eq:1050}
Y'(t) + \frac{1}{\We} \underbrace{\int_0^{+\infty} \!\! m(s) \partial_s \big( \|\bG\|_{L^2(\Omega)}^2 \big) \, \mathrm ds}_{\mathscr I}
\leq
\eps \, \|\nabla \bu\|_{L^2(\Omega)}^2 + C_Y(t) Y(t).
\end{equation}
Using a integration by part, the integral~$\mathscr I$ becomes
\begin{equation}
\mathscr I = \int_0^{+\infty}  \!\! -m'(s) \|\bG\|_{L^2(\Omega)}^2 \, \mathrm ds + \Big[ m(s) \|\bG\|_{L^2(\Omega)}^2 \Big]_0^{+\infty}.
\end{equation}
Using the following arguments:
\par\noindent
- The memory function $m$ is non-increasing, that is $-m'\geq 0$;\par\noindent
- The function~$\bG$ is bounded on $\R^+\! \times \! (0,T)\times \!\Omega$ and $m\in L^1(\R^+)$, that is $\dsp \lim_{s\to +\infty} m(s) \|\bG\|_{L^2(\Omega)}^2 = 0$;\par\noindent
- We have the following development with respect to the variable~$s$ for $\bG$: 
$$
\bG(s) = \bG\big|_{s=0} + s \, \partial_s \bG\big|_{s=0} + o(s) \sim s \, \We \, \nabla \bu.
$$
Moreover, $m\in L^1(\R^+)$ so that $\dsp \lim_{s\to 0} m(s) \|\bG\|_{L^2(\Omega)}^2 = 0$.\\[0.3cm]
Hence the integral~$\mathscr I$ is non-negative so that the estimate~\eqref{eq:1050} now reads
\begin{equation}\label{eq:unicity3}
Y'(t) \leq
\eps \, \|\nabla \bu\|_{L^2(\Omega)}^2 + C_Y(t) Y(t).
\end{equation}
{\bf Uniqueness result --}
Finally, adding \eqref{eq:unicity1} and~\eqref{eq:unicity3} with the choice $\eps=1-\omega$, and using the estimate~\eqref{eq:unicity2}, we obtain
\begin{equation}
(Y+Z)'(t) \leq C_{YZ}(t) \, (Y+Z)(t),
\end{equation}
where the function~$C_{YZ}$ is a linear combination of $C_Y$ and $C_Z$. In particular we have $C_{YZ}\in L^1(0,T)$.
The Gronwall lemma and the initial condition $Y(0)=Z(0)=0$ imply that $Y=Z=0$.
We deduce that $\bu=\b0$ and $\bG=\b0$ and that consequently the stress $\btau=\b0$ and the pressure~$p$ is constant in $(0,T)\times \Omega$.
\cqfd

\section{Proof for the global existence with small data}\label{part:global}

This section is devoted to the proof of the local existence theorem~\ref{th:global}.
Arguing as in the proof of the theorem~\ref{th:local}, we introduce the space~$\mathscr{B}(T)$, the subspaces $\mathscr{H}(T,R_1,R_2,R_3)$ and the mapping~$\Phi$.
\par\noindent
For $(\overline{\bu} , \overline{\bG}, \overline{\btau})\in \mathscr{H}(T,R_1,R_2,R_3)$ and $(\bu,\bG,\btau) = \Phi(\overline{\bu} , \overline{\bG} , \overline{\btau})$ recall that we have the following estimates (see the estimates~\eqref{eq:1499}, \eqref{eq:1112} and~\eqref{eq:1114}):
\begin{equation}\label{eq:1498}
\begin{aligned}
\|\bu\|_{L^r(0,T;W^{2,p}(\Omega))} + \| \partial_t \bu\|_{L^r(0,T;L^p(\Omega))}
\leq
\frac{C_1}{1-\omega} \Big( & \|\bu_0\|_{D^r_p(\Omega)} + \|\boldsymbol{f}\|_{L^r(0,T;L^p(\Omega))} \\
& + \Re \, C \, T^\alpha R_1^\beta \|\bu_0\|_{D^r_p(\Omega)}^\gamma + \Re \, C \, T^\delta R_1^2 + T^{\frac{1}{r}} R_3  \Big),
\end{aligned}
\end{equation}
\begin{equation}\label{eq:1112b}
\begin{aligned}
& \|\bG\|_{L^\infty(\R^+\times (0,T);W^{1,p}(\Omega))} + \|\partial_t \bG\|_{L^\infty(\R^+;L^r(0,T;L^p(\Omega)))} 
\leq
C_2 \big( 1+ R_1 \big) \mathrm{exp} \big( C_3 T^{1-\frac{1}{r}} R_1 \big).
\end{aligned}
\end{equation}
\begin{equation}\label{eq:1114b}
\begin{aligned}
\|\btau\|_{L^\infty(0,T;W^{1,p}(\Omega))} + \|\partial_t \btau\|_{L^r(0,T;L^p(\Omega))}
\leq
\frac{\omega}{\We} F_0(R_2).
\end{aligned}
\end{equation}
Note that the constants~$C_1$, $C_2$, $C_3$, $C$ and the function~$F_{0}$ introduced in these three estimates do not depend on~$\omega$.
For a time $T>0$ given, we successively choose
\begin{equation}
\begin{aligned}
& R_1^\star = \frac{1-\omega}{2C_1 \, \Re \, C\, T^\delta}, \\
& R_2^\star = C_2 \big( 1+ R_1^\star \big) \mathrm{exp} \big( C_3 T^{1-\frac{1}{r}} R_1^\star \big) + \|\bG_{\mathrm{old}}\|_{L^\infty(\R^+;W^{1,p}(\Omega))}, \\
& R_3^\star = \frac{\omega}{\We} F_0(R_2^\star).
\end{aligned}
\end{equation}
We verify that, for $\bu_0$, $\boldsymbol{f}$ small enough in their norms, and for $\omega$ small enough too, $\mathscr{H}(T,R_1^\star,R_2^\star,R_3^\star) \neq \emptyset$ (that is~\eqref{eq:2128} holds).
For such choices of $R_1^\star$, $R_2^\star$, $\omega$ and small norms of $\bu_0$, $\boldsymbol{f}$ and~$\bG_{\mathrm{old}}$, the inequalities~\eqref{eq:1498}, \eqref{eq:1112b} and~\eqref{eq:1114b} imply that $\Phi(\mathscr{H}(T,R_1^\star,R_2^\star,R_3^\star)) \subset \mathscr{H}(T,R_1^\star,R_2^\star,R_3^\star)$.
Moreover the function $\Phi$ is continuous and $\mathscr{H}(T,R_1^\star,R_2^\star,R_3^\star)$ is a convex compact subset of $\mathscr{B}(T)$.
We conclude the proof using the Schauder's theorem again.
\cqfd

\section{Proof on the stationary problem}\label{part:station}

The stationary problem associated to the system~\eqref{eq:1321} reads
\begin{equation}\label{eq:1320}
\left\{
\begin{aligned}
& \Re \, \bu\cdot \nabla \bu + \nabla p - (1-\omega) \Delta \bu = \div\, \btau + \boldsymbol f, \\
& \div\, \bu = 0, \\
& \btau(\bx) = \frac{\omega}{\We} \int_0^{+\infty} m(s) \, \mathcal S\big(\bG(s,\bx)\big)\, \mathrm ds,\\
& \frac{1}{\We} \partial_s \bG + \bu \cdot \nabla \bG = \bG \cdot \nabla \bu,
\end{aligned}
\right.
\end{equation}
coupled with the ``boundary'' conditions $\bu\big|_{\partial \Omega} = \b0$ and $\bG\big|_{s=0} = \bdelta$.\\[0.3cm]
We will introduce a Banach space $\mathscr{C}$, a convex compact subset $\mathscr{K}(R_1,R_3)$ and a continuous mapping $\Phi:\mathscr{K}(R_1,R_3) \longrightarrow \mathscr{C}$ in such a way that the system~\eqref{eq:1320} is equivalent to a fixed point equation for $\Phi$ in $\mathscr{K}(R_1,R_3)$.
More precisely, let us set $\mathscr{C}=W^{1,p}_0(\Omega) \times L^p(\Omega)$ and
\begin{equation}
\begin{aligned}
\mathscr{K}(R_1,R_3) = \Big\{ (\overline{\bu} , \overline{\btau})\in \mathscr{C} 
& ~ ; ~
\overline \bu \in D(A_p),
~~ \overline \btau \in W^{1,p}(\Omega), \\
& \|\overline \bu\|_{W^{2,p}(\Omega)} \leq R_1,
~~ \|\overline \btau\|_{W^{1,p}(\Omega)} \leq R_3 ~ \Big\}.
\end{aligned}
\end{equation}
\par\noindent
For each $R_1>0$ and $R_3>0$ it is obvious that we have $\mathscr{K}(R_1,R_3)\neq \emptyset$.
We now consider the application $\Phi : (\overline{\bu} , \overline{\btau}) \in \mathscr{K}(R_1,R_3) \longmapsto (\bu, \btau) \in \mathscr{B}$,
where $\bu$ is the unique solution of the stationary Stokes problem
\begin{equation}\label{eq:1441}
\left\{
\begin{aligned}
& \nabla p - (1-\omega) \Delta \bu = \overline \bg, \\
& \div\, \bu = 0, \\
& \bu\big|_{\partial \Omega} = \b0,
\end{aligned}
\right.
\end{equation}
with the following source term $\overline \bg = -\Re\, \overline{\bu} \cdot \nabla \overline{\bu} + \div\, \overline{\btau} + \boldsymbol{f}$;
and where $\btau$ is given by the integral formula
\begin{equation}\label{eq:2056}
\btau(\bx) = \frac{\omega}{\We} \int_0^{+\infty} m(s) \, \mathcal S\big(\bG(s,\bx)\big)\, \mathrm ds,
\end{equation}
the tensor~$\bG$ being the unique solution of the equation
\begin{equation}\label{eq:2055}
\left\{
\begin{aligned}
& \frac{1}{\We} \partial_s \bG + \overline \bu \cdot \nabla \bG = \bG \cdot \nabla \overline \bu, \\
& \bG\big|_{s=0} = \bdelta.
\end{aligned}
\right.
\end{equation}

\begin{remark}
In the proof of the theorem~\ref{th:station}, for reasons of clarity, we often only indicate the dependence of the constants with respect to the parameter~$\omega$.
We write $A\lesssim B$ if there exists a constant $C$ which does not depend on~$\omega$ (but depending possibly on $\Omega$, $p$, $d$,...) such that $A\leq CB$.
\end{remark}

\subsection{Estimates for the velocity field~$\bu$}

It is well known (the proof is given in~\cite{Ladyzhen}) that the unique solution of the Stokes problem~\eqref{eq:1441} belongs in $W^{2,p}(\Omega)$ as soon as the source term $\overline \bg\in L^p(\Omega)$.
Moreover we have the following continuity estimate:
\begin{equation}
\|\nabla \bu\|_{W^{1,p}(\Omega)}\leq \frac{C_4}{1-\omega} \|\overline \bg\|_{L^p(\Omega)},
\end{equation}
where the constant $C_4$ only depends on the domain~$\Omega$ and the integer~$p$.
For our purposes, the expression of the source term~$\overline \bg$ allows us to deduce
\begin{equation}\label{eq:1609}
\|\nabla \bu\|_{W^{1,p}(\Omega)}\leq \frac{C_4}{1-\omega} \Big( \|\boldsymbol{f}\|_{L^p(\Omega)} + \Re \|\nabla \overline \bu\|_{W^{1,p}(\Omega)}^2 + \|\overline \btau\|_{W^{1,p}(\Omega)} \Big).
\end{equation}

\subsection{Estimates for the deformation gradient~$\bG$}

Following the same method that those to obtain an estimate of~$\bG$ for the unsteady case (see the subsection~\ref{part:estimates-1}), we prove that the solution~$\bG$ of the system~\eqref{eq:2055} (which exists and is unique by the characteristic method) satisfies
\begin{equation}
\frac{1}{\We} d_s \big( \|\bG\|_{W^{1,p}(\Omega)} \big)
\leq 3C_0 \|\nabla \overline \bu\|_{W^{1,p}(\Omega)} \|\bG\|_{W^{1,p}(\Omega)},
\end{equation}
where $C_0=C_0(\Omega,p)$ corresponds to a constant of the continuous injection $W^{1,p}(\Omega) \hookrightarrow L^\infty(\Omega)$.
Since $\|\bG\|_{W^{1,p}(\Omega)}|_{s=0}=\sqrt{d}\, |\Omega|^{\frac{1}{p}}$, the classical Gronwall lemma implies that for all $s\geq 0$,
\begin{equation}\label{eq:1607}
\|\bG(s,\cdot)\|_{W^{1,p}(\Omega)} \leq \sqrt{d}\, |\Omega|^{\frac{1}{p}} \mathrm{exp} \Big( 3C_0 \We \|\nabla \overline \bu\|_{W^{1,p}(\Omega)} s \Big).
\end{equation}
To give sense to the ``initial'' condition $\bG|_{s=0}=\bdelta$ we need estimate on~$\partial_s \bG$.
Isolating this term on the equation~\eqref{eq:2055} we obtain
\begin{equation}
\|\partial_s \bG(s,\cdot)\|_{L^p(\Omega)} \leq \We \Big( \|\bG(s,\cdot)\|_{W^{1,p}(\Omega)} \|\nabla \overline \bu\|_{L^p(\Omega)} + \|\overline \bu\|_{W^{1,p}(\Omega)} \|\nabla \bG(s,\cdot)\|_{L^p(\Omega)}\Big).
\end{equation}
That is, due to the estimate~\eqref{eq:1607}, for all $s\in \R^+$:
\begin{equation}\label{eq:2055b}
\|\partial_s \bG(s,\cdot)\|_{L^p(\Omega)} \lesssim \|\nabla \overline \bu\|_{L^p(\Omega)} \mathrm{exp} \Big( 3C_0 \We \|\nabla \overline \bu\|_{W^{1,p}(\Omega)} s \Big).
\end{equation}

\subsection{Estimates for the stress tensor~$\btau$}

To control the stress tensor given by the integral~\eqref{eq:2056}, we must control the quantity~$\mathcal S(\bG)$.
Using the assumptions (A4) on the growth of this function: $|\mathcal S(\bG)|\lesssim |\bG|^a$, $|\mathcal S'(\bG)|\lesssim |\bG|^b$ we deduce that for all $\bx\in \Omega$:
\begin{equation}
|\btau(\bx) |^p
\lesssim \omega^p \! \int_0^\infty m(s) \, |\bG(s,\bx)|^{ap}\, \mathrm ds,
\end{equation}
and
\begin{equation}
|\nabla \btau(\bx) |^p
\lesssim \omega^p \! \int_0^\infty m(s) \, |\nabla \bG(s,\bx)|^p \, |\bG(s,\bx)|^{bp}\, \mathrm ds.
\end{equation}
Since $\bG \in W^{1,p}(\Omega) \subset L^\infty(\Omega)$ (with respect to the spatial variable), we have
\begin{equation}
\|\btau\|_{L^p(\Omega)}^p
\lesssim \omega^p \! \int_0^\infty m(s) \, \|\bG(s,\cdot)\|_{W^{1,p}(\Omega)}^{ap}\, \mathrm ds,
\end{equation}
and
\begin{equation}
\|\nabla \btau\|_{L^p(\Omega)}^p
\lesssim \omega^p \! \int_0^\infty m(s) \, \|\bG(s,\cdot)\|_{W^{1,p}(\Omega)}^{(b+1)p}\, \mathrm ds.
\end{equation}
Using the estimate~\eqref{eq:1607} on~$\bG$, and the assumption (A3), that is $m(s)\lesssim \mathrm e^{-\alpha s}$, we obtain
\begin{equation}\label{eq:1608}
\|\btau \|_{W^{1,p}(\Omega)}
\leq \omega \, C_5 \int_0^{\infty} \mathrm{exp}\big( (3C_0 \We \, c \, p\, \|\overline \bu\|_{W^{1,p}(\Omega)} - \alpha) s \big) \, \mathrm ds,
\end{equation}
where $c=\max\{a,b+1\}$ and where $C_5$ does not depend on~$\omega$.

\subsection{Proof of the existence result for the theorem~\ref{th:station}}

Let $(\overline{\bu} , \overline{\btau})\in \mathscr{K}(R_1,R_3)$ and $(\bu,\btau)=\Phi(\overline{\bu} , \overline{\btau})$.
Using the relations~\eqref{eq:1609} and~\eqref{eq:1608}, we have
\begin{equation}\label{eq:1609b}
\|\nabla \bu\|_{W^{1,p}(\Omega)}\leq \frac{C_4}{1-\omega} \Big( \|\boldsymbol{f}\|_{L^p(\Omega)} + \Re \, R_1^2 + R_3 \Big),
\end{equation}
\begin{equation}\label{eq:1608b}
\|\btau \|_{W^{1,p}(\Omega)}
\leq \omega \, C_5 \int_0^{\infty} \mathrm e^{(3C_0 \We \, c \, p \, R_1 - \alpha) s } \, \mathrm ds.
\end{equation}
Taking $R_1$ small enough to have $3C_0 \We \, c \, p \, R_1 < \alpha/2$, the inequality~\eqref{eq:1608b} becomes
\begin{equation}\label{eq:1608bb}
\|\btau \|_{W^{1,p}(\Omega)}
\leq \frac{2\omega C_5}{\alpha}.
\end{equation}
We now choose $R_3=\frac{2 \omega C_5}{\alpha}$ and we verify that for~$R_1$, $\boldsymbol{f}$ and~$\omega$ small enough, the relations~\eqref{eq:1609b} and~\eqref{eq:1608bb} imply the inclusion~$\Phi(\mathscr{K}(R_1,R_3))\subset \mathscr{K}(R_1,R_3)$.
As for the proof of the theorem~\ref{th:local}, the Schauder'theorem implies the existence of a solution with estimates for~$\bu$ and for~$\btau$.
Next, using~\eqref{eq:1607} and~\eqref{eq:2055b}, we deduce the estimates for~$\bG$ and for~$\partial_s \bG$.

\subsection{Proof of the uniqueness result for the theorem~\ref{th:station}}

Proceeding as in the proof of the theorem~\ref{th:unic} (see section~\ref{part:unic}), we consider two solutions $(\bu_i,\btau_i) \in \mathscr{K}(R_1,R_3)$ with $i\in\{1,2\}$.
We note $\bG_i$ the corresponding deformation gradients.
From the estimate~\eqref{eq:1607}, we have for all $s\in \R^+$:
\begin{equation}\label{eq:1607b}
\|\bG_i(s,\cdot)\|_{W^{1,p}(\Omega)} \leq \sqrt{d}\, |\Omega|^{\frac{1}{p}} \mathrm e^{3C_0 \We R_1 s}.
\end{equation}
We next introduce the differences $\bu=\bu_1-\bu_2$, $\btau = \btau_1-\btau_2$ and $\bG=\bG_1-\bG_2$ and consider the equations satisfied by these differences.
\par\noindent
$\bullet$
The scalar product of the Stokes equation for~$\bu$ by~$\bu$ in~$L^2(\Omega)$ gives
\begin{equation}\label{eq:unicity1station}
(1-\omega) \|\nabla \bu\|_{L^2(\Omega)}^2
\leq \frac{4}{1-\omega}\|\btau\|_{L^2(\Omega)}^2 + \Re \, \|\nabla \bu_2\|_{L^\infty(\Omega)} \|\bu\|_{L^2(\Omega)}^2.
\end{equation}
\par\noindent
$\bullet$
From the definition of~$\btau$ as a difference of two integrals depending on $\mathcal S(\bG_1)$ and $\mathcal S(\bG_2)$ we obtain
\begin{equation}
|\btau(\bx)| \leq \frac{\omega}{\We} \int_0^{\infty} m(s) \, \big| f'(\bG_3(s,\bx)) \big| \, |\bG(s,\bx)|\, \mathrm ds,
\end{equation}
where $\dsp \|\bG_3(s,\cdot)\|_{W^{1,p}(\Omega)} \leq \max_{i\in \{1,2\}}\{\|\bG_i(s,\cdot)\|_{W^{1,p}(\Omega)}\} \lesssim \mathrm e^{3C_0 \We R_1 s}$.
Using the bound $|f'(\bG_3)| \lesssim |\bG_3|^b$ we obtain
\begin{equation}
|\btau(\bx)| \lesssim \omega \int_0^{\infty} m(s) |\bG_3(s,\bx)|^b \, |\bG(s,\bx)|\, \mathrm ds.
\end{equation}
By the continuous injection $W^{1,p}(\Omega) \hookrightarrow L^\infty(\Omega)$, and from the estimate on $\|\bG_i\|_{W^{1,p}(\Omega)}$, $i\in \{1,2,3\}$, we have
\begin{equation}
|\btau(\bx)| \lesssim \omega \int_0^{\infty} m(s) \, \mathrm e^{3C_0 \We R_1 b s} \, |\bG(s,\bx)|\, \mathrm ds.
\end{equation}
Taking the square and integrating with respect to $\bx$ (after using the Hölder inequality) we deduce
\begin{equation}\label{eq:unicity2station}
\|\btau\|_{L^2(\Omega)}^2 \lesssim \omega^2 \int_0^{\infty} m(s) \, \mathrm e^{6C_0 \We R_1 b s}\, \|\bG(s,\cdot)\|_{L^2(\Omega)}^2 \, \mathrm ds.
\end{equation}
\par\noindent
$\bullet$
Taking the scalar product of the equation on~$\bG$ by~$\bG$ in $L^2(\Omega)$, we obtain (see~\eqref{eq:2205} for similar calculations)
\begin{equation}
\begin{aligned}
\frac{1}{2 \We} d_s \big( \|\bG\|_{L^2(\Omega)}^2 \big)
\leq
& \|\bG_1\|_{L^\infty(\Omega)} \|\nabla \bu\|_{L^2(\Omega)} \|\bG\|_{L^2(\Omega)} \\
& \qquad + \|\nabla \bu_2\|_{L^\infty(\Omega)} \|\bG\|_{L^2(\Omega)}^2
+ \|\nabla \bG_2\|_{L^d(\Omega)} \|\nabla \bu\|_{L^2(\Omega)} \|\bG\|_{L^2(\Omega)}.
\end{aligned}
\end{equation}
From the Young inequality, this relation is written
\begin{equation}
d_s \big( \|\bG\|_{L^2(\Omega)}^2 \big)
\leq
A(s) \|\nabla \bu\|_{L^2(\Omega)}^2 + B \|\bG\|_{L^2(\Omega)}^2,
\end{equation}
where $A(s) =\frac{\We}{R_1} \big( \|\bG_1\|_{L^\infty(\Omega)}^2 + \|\nabla \bG_2\|_{L^d(\Omega)}^2 \big)$ and $B= 2 \, \We \big( \|\nabla \bu_2\|_{L^\infty(\Omega)} +  R_1 \big)$.
By the Gronwall lemma (and using the fact that $\bG(0,\bx)=\b0$) we obtain
\begin{equation}\label{eq:2320?}
\|\bG(s,\cdot)\|_{L^2(\Omega)}^2
\leq \mathrm e^{Bs}
\int_0^s \mathrm e^{-Bs'} A(s') \|\nabla \bu\|_{L^2(\Omega)}^2 \, \mathrm ds'.
\end{equation}
Using the bound on $\bG_i$ we have $A(s')\lesssim \frac{1}{R_1} \mathrm e^{6C_0 \We R_1 s}$ for all $s'\in \R^+$.
We deduce
\begin{equation}
\|\bG(s,\cdot)\|_{L^2(\Omega)}^2
\lesssim \frac{\mathrm e^{Bs}}{R_1} 
\Big( \int_0^s \mathrm e^{(6C_0 \We R_1-B)s'} \mathrm ds' \Big) \|\nabla \bu\|_{L^2(\Omega)}^2.
\end{equation}
Now, the estimate~\eqref{eq:unicity2station} becomes
\begin{equation}
\|\btau\|_{L^2(\Omega)}^2 \lesssim \frac{\omega^2}{R_1} \Big( \underbrace{\int_0^{\infty} \int_0^s m(s) \, \mathrm e^{(6C_0 \We R_1 b +B)s}\, \mathrm e^{(6C_0 \We R_1-B)s'} \mathrm ds' \, \mathrm ds}_{\mathscr I} \Big) \|\nabla \bu\|_{L^2(\Omega)}^2.
\end{equation}
Assuming (A3), that is $m(s)\lesssim \mathrm e^{-\alpha s}$, the integral $\mathscr I$ satisfies
\begin{equation}
\mathscr I \lesssim \int_0^\infty \mathrm e^{(6C_0 \We R_1-B)s'} \Big( \int_{s'}^\infty \mathrm e^{(6C_0 \We R_1 b +B - \alpha)s} \mathrm ds \Big) \mathrm ds'.
\end{equation}
Hence, for $R_1$ small enough (such that $2\We (3C_0 b + 2) R_1 < \alpha$, we will note that $B\leq 4\We \, R_1$) the integration with respect to the variable $s$ in the last integral converges and we have
\begin{equation}
\mathscr I \lesssim \frac{1}{\alpha - (6C_0 \We R_1 b + B)} \int_0^\infty \mathrm e^{(6C_0 \We R_1(b+1)-\alpha)s'} \mathrm ds'.
\end{equation}
Finally, if in addition we choose~$R_1$ small enough to have $6C_0 \We R_1(b+1)<\alpha$ then the integral~$\mathscr I$ converges and we have the following bound for the stress:
\begin{equation}\label{eq:1425}
\|\btau\|_{L^2(\Omega)}^2 \lesssim \frac{\omega^2}{R_1} \|\nabla \bu\|_{L^2(\Omega)}^2.
\end{equation}
\par\noindent
Using the relations~\eqref{eq:unicity1station} and~\eqref{eq:1425}, we deduce that for $R_1$ and~$\omega$ small enough we have $\bu=\b0$ and $\btau=\b0$.
Next, the relation~\eqref{eq:2320?} implies $\bG=\b0$, that concludes the uniqueness proof of the theorem~\ref{th:station}.
\cqfd

\section{Conclusion}\label{part:conclusion}

In this article we are interested in the mathematical properties of models of viscoelastic flow.
We have shown that many known results for differential laws could be adapted to integral models.
Nevertheless some differences persist:
\begin{itemize}
\item[$\bullet$] Our results are formulated in the $L^r\!-\!L^p$ context (following~\cite{Fernandez-Guillen-Ortega}). It seems possible to reformulate them for more regular solutions in a $H^s$ context (following~\cite{Guillope-Saut3}).
\item[$\bullet$] Some results: the global existence with small data (theorem~\ref{th:global}) and the stationary study (theorem~\ref{th:station}), are proved for the relaxation parameter~$\omega$ small enough only, that is to say for the flows which are not too elastic.
About the differential models, this assumption can be removed, see for instance~\cite{Chupin2, Molinet} for the global existence, and~\cite{Fernandez-Guillen-Ortega} for the stationary problem.
In these two cases, the results about differential models strongly use the structure of the equation, and it seems difficult to adapt such methods for integral models (see also the remark~\ref{rem:theorem}).
\item[$\bullet$] There exists differential models which have no apparent equivalent in terms of integral models, for instance the co-rotational Oldroyd model.
This study does not cover these cases.
Similarly, there are also integral models more general than those studied here.
In~\cite{Tanner}, R.~I. Tanner introduce models where the memory~$m$ also depends on the invariants~$I_1$ and~$I_2$.
It might be interesting to study from a theoretical point of view these models and to observe if the approach taken here fits well.
\end{itemize}
On the other hand, it is possible that classical integral models perform better than differential models of Oldroyd type.
In fact, most of these models have a stress which is naturally bounded {\it via} the definition of $\mathcal S$, see the examples given by the equation~\eqref{ex:0001}, \eqref{ex:0002} or~\eqref{ex:0003}, and the Appendix~\ref{part:appendix0}.
While getting weak solution seems very difficult, knowing {\it a priori} bound on the stress is an interesting information (see~\cite{Chemin} for an example of a criteria for the explosion in the Oldroyd model).\\[0.3cm]
Finally, the theoretical results shown in this paper allow to consider a lot of work on models of integral type.
The well-known results for the differential model can be generalized to integral models.
For instance, the one dimensional shearing motions and Poiseuille flows admit global existence for usual differential models, see~\cite{Guillope-Saut4}.
In this regards, the works of A.C.T.~Aarts and A.A.F.~van de Ven~\cite{Aarts} are interesting: they study the Poiseuille flow of a K-BKZ model.
Would it be possible to prove global existence for such one dimensional flows when we use more general integral models\,?
An other possible generalization concerns the behavior of viscoelastic flows in thin geometries (in the fields such as polymer extrusion or lubrication), or in thin free-surface flows (for study mudslide or oil slick).
Recent works~\cite{Bayada-Chupin-Grec-viscoelastique,Bayada-Chupin-Martin} and~\cite{Bouchut} can provide answers to the differential models, and we can imagine the same kind of works for integral models.

\appendix

\section{Tensors and the strain measure function~$\mathcal S$}\label{part:appendix0}

\subsection{Some remarks on the invariants of the Finger tensor ($d=3$)}

For a matrix $\bB\in \mathcal L(\R^3)$, we usually define three invariants:
$$
I_1=\mathrm{Tr}(\bB),
\quad
I_2=\frac{1}{2}\big( (\mathrm{Tr}(\bB))^2 - \mathrm{Tr}(\bB^2) \big),
\quad
I_3=\det \bB.
$$
We specify in this subsection some properties of these invariants in the context studied here, that is when~$\bB$ represents a Finger tensor of an incompressible flow.
\par\noindent
First of all, this incompressibility condition implies that $\det \bB=1$.
Consequently the third invariant $I_3$ is useless.
Next, using the Cayley-Hamilton theorem we have $\bB^{-1} = \bB^2-I_1\bB+I_2\, \bdelta$ and we deduce that
$$I_2=\mathrm{Tr}(\bB^{-1}).$$
By definition, $\bB=\transp{\bF}\cdot \bF$ is real positive-definite matrix, and consequently it is diagonalizable.
Using a basis formed by its eigenvectors, we have
$
I_1=\lambda_1 + \lambda_2 + \lambda_3
$
whereas $\lambda_1\lambda_2\lambda_3=1$.
An inequality of arithmetic and geometric means indicates that $I_1\geq 3$, and in a similar way we prove that $I_2\geq 3$.
We deduce that for all $\beta\in [0,1]$ we have 
$$\beta I_1 + (1-\beta) I_2 \geq 3.$$
That mathematically justifies the PSM model~\eqref{ex:0001} and the Wagner model~\eqref{ex:0002}.

\subsection{Notion of derivative for the Strain measure tensor}

The strain measure function~$\mathcal S : \mathcal L(\R^d) \longrightarrow \mathcal L(\R^d)$ can be viewed as an application $\R^{d\times d} \longrightarrow \R^{d\times d}$.
For $1\leq i,j\leq d$, we introduce the matrix
$$\bE_{ij} = 
\begin{array}{cl}
\begin{pmatrix} 
~~0 & & 0~ \\[0.1cm]
 &1&\cdots \\[-0.1cm]
~~0 & \vdots & 0~ \\
\end{pmatrix} & \hspace{-0.3cm}j\\
i & \\
\end{array}$$
as element of the basis of $\mathcal L(\R^d)$ and we define the differential of $\mathcal S$ from its Jacobian $\dsp \mathcal S' := \Big( \frac{\partial \mathcal S_{k\ell}}{\partial \bE_{ij}} \Big)_{ijk\ell}$.
Note that in this paper, we use the notation $\dsp \partial_{(ij)} \mathcal S_{k\ell} := \frac{\partial \mathcal S_{k\ell}}{\partial \bE_{ij}}$.

\subsection{Norm for the $4$-tensor}

The notion of derivative introduced above involves the use of tensors of order~$4$.
Recall that for a $2$-tensor $\bA=(\bA)_{ij}$, we use the usual algebra norm defined by $|\bG|^2:=\mathrm{Tr}(\transp{\bG} \cdot \bG)$.
For a $4$-tensors $\bH=(\bH)_{ijk\ell}$, we introduce the following algebra norm:
$$|\bH|^2:=\sum_{i,j,k,\ell}\bH_{ijk\ell}^2.$$
We will note that this norm having the following property: $|\bA \otimes \bB| = |\bA|\, |\bB|$ for any $2$-tensors $\bA$ and $\bB$.

\subsection{Example for a PSM model}

Consider the example function corresponding to a PSM model (see the equation~\eqref{ex:0001} with $\alpha=4$ and $\beta=1$):
\begin{equation}
\mathcal S:\bG \in \mathcal L(\R^d) \longmapsto \frac{\bB}{1+\mathrm{Tr}(\bB)} \in \mathcal L(\R^d)
\qquad \text{where}\quad \bB = \transp{\bG} \cdot \bG.
\end{equation}
\begin{proposition}
For all $\bG\in \mathcal L(\R^d)$ we have $\dsp |\mathcal S(\bG)|\leq 1$.
\end{proposition}
\proof
Using the norm on the $2$-tensor, we have for all $\bG \in \mathcal L(\R^d)$:
\begin{equation}
|\mathcal S(\bG)| = \frac{|\transp{\bG} \cdot \bG|}{1+|\bG|^2} \leq \frac{|\bG|^2}{1+|\bG|^2},
\end{equation}
that implies that $\mathcal S$ is bounded by the constant~$1$.
\cqfd
\begin{proposition}
For all $\bG\in \mathcal L(\R^d)$ we have $\dsp |\mathcal S'(\bG)|\leq \frac{2(1+\sqrt{d})}{|\bG|}$.
\end{proposition}
\proof
The derivative of the function $\mathcal S$ is a function with values in the set of $4$-tensors:
\begin{equation}
\begin{aligned}
\mathcal S'(\bG)_{ijk\ell} 
& = \partial_{(ij)}\Big( \frac{\bB_{k\ell}}{1+\mathrm{Tr}(\bB)} \Big) \\
& = \frac{\partial_{(ij)}(\bB_{k\ell})}{1+\mathrm{Tr}(\bB)} - \frac{\partial_{(ij)}(\mathrm{Tr}(\bB)) \bB_{kl}}{(1+\mathrm{Tr}(\bB))^2} \\
& = \frac{(\bdelta_{kj}\bG_{i\ell} + \bdelta_{\ell j}\bG_{ik})}{1+\mathrm{Tr}(\bB)} - \frac{2\bG_{ij} \bB_{kl}}{(1+\mathrm{Tr}(\bB))^2},
\end{aligned}
\end{equation}
Taking the $4$-tensor norm, we deduce that
\begin{equation}
|\mathcal S'(\bG)| \leq \frac{2|\bdelta|\, |\bG|}{1+|\bG|^2} + \frac{2|\bG|^3}{(1+|\bG|^2)^2}
\end{equation}
Since $|\bdelta|=\sqrt{d}$, this implies that $|\bG|\, |\mathcal S'(\bG)|$ is bounded by $2(1+\sqrt{d})$.
\cqfd

\section{Gronwall type lemma}\label{part:appendix1}

\begin{lemma}\label{lem:gronwall}
Let $f:\R^+ \mapsto \R^+$ a positive and locally integrable function.
If a function $y:\R^+ \! \times \R^+ \mapsto \R$ satisfies, for all $(s,t)\in \R^+ \! \times \R^+$:
\begin{equation}\label{eq:lem1}
\partial_t y(s,t) + \frac{1}{\We} \partial_s y(s,t) \leq f(t)\, y(s,t)
\end{equation}
then we have, for all $(s,t)\in \R^+ \! \times \R^+$:
\begin{equation}\label{eq:lem2}
y(s,t) \leq  \zeta(s,t) \, \mathrm{exp}\bigg( \int_0^t f(t')\, \mathrm dt' \bigg),
\end{equation}
where $\dsp \zeta(s,t) = \left\{
\begin{aligned}
& y\big( s-\frac{t}{\We},0\big) \quad \text{if $t\leq \We\, s$}, \\
& y(0,t-\We\, s) \quad \text{if $t > \We\, s$}.
\end{aligned}
\right.
$
\end{lemma}

\proof
Introducing the new variables $u=\frac{1}{2}(\We\, s+t)$ and $v=\frac{1}{2}(\We\, s-t)$, we can write the first equation of the system~\eqref{eq:lem1} as a system on the function $z(u,v)=y(s,t)$:
\begin{equation}
\partial_u z(u,v) \leq f(u-v) \, z(u,v).
\end{equation}
Since the function $f$ is locally integrable, we obtain
\begin{equation}
\partial_u \bigg[ z(u,v)\, \mathrm{exp}\bigg( - \int_0^{u-v} f(t')\, \mathrm dt' \bigg) \bigg] \leq 0.
\end{equation}
Integrating this relation between $|v|$ and $u$, we deduce
\begin{equation}
z(u,v)\, \mathrm{exp}\bigg( - \int_0^{u-v} f(t')\, \mathrm dt' \bigg) \leq z(|v|,v)\, \mathrm{exp}\bigg( - \int_0^{|v|-v} f(t')\, \mathrm dt' \bigg).
\end{equation}
Due to the positivity of the function~$f$, the exponential term in the last equation being less than~$1$.
According to the sign of $v$, we have $z(|v|,v)=y(0,t-\We\, s)$ or $z(|v|,v)=y(s-\frac{t}{\We},0)$.
That implies the result~\eqref{eq:lem2} announced in the lemma.
\cqfd


\noindent
{\bf Acknowledgments} - The author has been partially supported by the ANR project ANR-08-JCJC-0104-01 : RUGO (Analyse et calcul des effets de rugosit\'es sur les \'ecoulements).


\nocite{*}
\bibliographystyle{plain}
\bibliography{biblio}

\begin{thebibliography}{10}

\bibitem{Aarts}
A.~C.~T. Aarts and A.~A.~F. van~de Ven.
\newblock Transient behaviour and stability points of the {P}oiseuille flow of
  a {KBKZ}-fluid.
\newblock {\em J. Engrg. Math.}, 29(4):371--392, 1995.

\bibitem{Bayada-Chupin-Martin}
G.~Bayada, L.~Chupin, and S.~Martin.
\newblock Viscoelastic fluids in a thin domain.
\newblock {\em Quart. Appl. Math.}, 65(4):625--651, 2007.

\bibitem{Bayada-Chupin-Grec-viscoelastique}
Guy Bayada, Laurent Chupin, and B{\'e}r{\'e}nice Grec.
\newblock Viscoelastic fluids in thin domains: a mathematical proof.
\newblock {\em Asymptot. Anal.}, 64(3-4):185--211, 2009.

\bibitem{KBKZ1}
B.~Bernstein, E.A. Kearsley, and L.J. Zapas.
\newblock A study of stress relaxation with finite strain.
\newblock {\em Trans. Soc. Rheol.}, 12:623--727, 1959.

\bibitem{KBKZ2}
B.~Bernstein, E.A. Kearsley, and L.J. Zapas.
\newblock A study of stress relaxation with finite strain.
\newblock {\em Trans. Soc. Rheol.}, 17:35--92, 1964.

\bibitem{Bird}
R.B. Bird, O.~Hassager, R.C. Armonstrong, and C.F. Curtiss.
\newblock {\em Dynamics of Polymeric Fluids}, volume~2 of {\em Kinetic Theory}.
\newblock John Wiley and Sons, New York, 1977.

\bibitem{Bouchut}
F.~Bouchut and S.~Boyaval.
\newblock A new model for shallow elastic fluids.
\newblock {\em preprint}, 2011.

\bibitem{Chambon}
F.~Chambon and H.H. Winter.
\newblock {\em J. Rheol.}, 31:683, 1987.

\bibitem{Chemin}
Jean-Yves Chemin and Nader Masmoudi.
\newblock About lifespan of regular solutions of equations related to
  viscoelastic fluids.
\newblock {\em SIAM J. Math. Anal.}, 33(1):84--112 (electronic), 2001.

\bibitem{Chupin2}
L.~Chupin.
\newblock Some theoretical results concerning diphasic viscoelastic flows of
  the {O}ldroyd kind.
\newblock {\em Adv. Differential Equations}, 9(9-10):1039--1078, 2004.

\bibitem{Coleman}
B.D. Coleman and V.J. Mizel.
\newblock {\em Arch Ration Mech Anal}, 29:18--31, 1968.

\bibitem{Currie}
P.K. Currie.
\newblock Constitutive equations for polymer melts predicted by the doi-edwards
  and curtiss-bird kinetic theory models.
\newblock {\em J. Non-Newtonian Fluid Mech.}, 11:53--68, 1982.

\bibitem{Doi}
M.~Doi and S.F. Edwards.
\newblock {\em The theory of polymer dynamics}.
\newblock Oxford University Press, 1988.

\bibitem{Fernandez-Guillen-Ortega-CRAS}
E.~Fern{\'a}ndez-Cara, F.~Guill{\'e}n, and R.R. Ortega.
\newblock Existence et unicit\'e de solution forte locale en temps pour des
  fluides non newtoniens de type {O}ldroyd (version {$L\sp s$}--{$L\sp r$}).
\newblock {\em C. R. Acad. Sci. Paris S\'er. I Math.}, 319(4):411--416, 1994.

\bibitem{Fernandez-Guillen-Ortega}
E.~Fern{\'a}ndez-Cara, F.~Guill{\'e}n, and R.R. Ortega.
\newblock Some theoretical results concerning non-{N}ewtonian fluids of the
  {O}ldroyd kind.
\newblock {\em Ann. Scuola Norm. Sup. Pisa Cl. Sci. (4)}, 26(1):1--29, 1998.

\bibitem{Freed-Diethelm}
A.D. Freed and K.~Diethelm.
\newblock Fractional calculus in biomechanics: a 3d viscoelastic model using
  regularized fractional derivative kernels with application to the human
  calcaneal fat pad.
\newblock {\em Biomechan Model Mechanobiol}, 5:203--215, 2006.

\bibitem{Friedman}
A.~Friedman.
\newblock Partial differential equations.
\newblock {\em Holt-Rinehart-Winston, New York}, 1976.

\bibitem{Giga}
Yoshikazu Giga and Hermann Sohr.
\newblock Abstract {$L^p$} estimates for the {C}auchy problem with applications
  to the {N}avier-{S}tokes equations in exterior domains.
\newblock {\em J. Funct. Anal.}, 102(1):72--94, 1991.

\bibitem{Guillope-Saut4}
C.~Guillop{\'e} and J.-C. Saut.
\newblock Global existence and one-dimensional nonlinear stability of shearing
  motions of viscoelastic fluids of {O}ldroyd type.
\newblock {\em RAIRO Mod\'el. Math. Anal. Num\'er.}, 24(3):369--401, 1990.

\bibitem{Guillope-Saut-CRAS}
C.~Guillop{\'e} and J.C. Saut.
\newblock R\'esultats d'existence pour des fluides visco\'elastiques \`a loi de
  comportement de type diff\'erentiel.
\newblock {\em C. R. Acad. Sci. Paris S\'er. I Math.}, 305(11):489--492, 1987.

\bibitem{Guillope-Saut3}
C.~Guillop{\'e} and J.C. Saut.
\newblock Existence results for the flow of viscoelastic fluids with a
  differential constitutive law.
\newblock {\em Nonlinear Anal.}, 15(9):849--869, 1990.

\bibitem{Guillope-Saut1}
C.~Guillop{\'e} and J.C. Saut.
\newblock Mathematical problems arising in differential models for viscoelastic
  fluids.
\newblock In {\em Mathematical topics in fluid mechanics (Lisbon, 1991)},
  volume 274 of {\em Pitman Res. Notes Math. Ser.}, pages 64--92. Longman Sci.
  Tech., Harlow, 1992.

\bibitem{Jonscher}
A.K. Jonscher.
\newblock {\em Nature}, 267:673, 1977.

\bibitem{Kaye}
A.~Kaye.
\newblock Non-newtonian flowin incompressible fluids.
\newblock {\em Technical Report}, 134, 1962.

\bibitem{Keunings}
R.~Keunings.
\newblock Finite element methods for integral viscoelastic fluids.
\newblock {\em Rheology Reviews}, pages 167--195, 2003.

\bibitem{Kim}
Jong~Uhn Kim.
\newblock Global smooth solutions of the equations of motion of a nonlinear
  fluid with fading memory.
\newblock {\em Arch. Rational Mech. Anal.}, 79(2):97--130, 1982.

\bibitem{Kjartansson}
E.~Kjartansson.
\newblock {\em J. Geophys. Res.}, 84:4737, 1979.

\bibitem{Ladyzhen}
O.A. Ladyzhenskaya.
\newblock The mathematical theory of viscous incompressible flow.
\newblock {\em Gordon and Breach, New York}, 1969.

\bibitem{Lions-Masmoudi-viscoelastique}
P.L. Lions and N.~Masmoudi.
\newblock Global solutions for some {O}ldroyd models of non-{N}ewtonian flows.
\newblock {\em Chinese Ann. Math. Ser. B}, 21(2):131--146, 2000.

\bibitem{Molinet}
L.~Molinet and R.~Talhouk.
\newblock On the global and periodic regular flows of viscoelastic fluids with
  a differential constitutive law.
\newblock {\em NoDEA Nonlinear Differential Equations Appl.}, 11(3):349--359,
  2004.

\bibitem{Oldroyd}
J.G. Oldroyd.
\newblock On the formulation of rheological equations of state.
\newblock {\em Proc. Roy. Soc. London. Ser. A.}, 200:523--541, 1950.

\bibitem{Papanastasiou}
A.C. Papanastasiou, L.~Scriven, and C.~Macosko.
\newblock {\em Journal of Rheology}, 27:381--410, 1983.

\bibitem{Renardy}
M.~Renardy.
\newblock An existence theorem for model equations resulting from kinetic
  theories of polymer solutions.
\newblock {\em SIAM J. Math. Anal.}, 22(2):313--327, 1991.

\bibitem{Scher}
H.~Scher and E.W. Montroll.
\newblock {\em Phys. Rev. B}, 12:2455, 1975.

\bibitem{Schiessel}
H.~Schiessel and A.~Blumen.
\newblock {\em Macromolecules}, 28:4013, 1995.

\bibitem{Tanner}
R.I. Tanner.
\newblock {\em Engineering Rheology}.
\newblock Clarendon Press, 1988.

\bibitem{Wagner0}
M.H. Wagner.
\newblock {\em Rheologica Acta}, 15:136, 1976.

\bibitem{Wagner1}
M.H. Wagner.
\newblock {\em Rheologica Acta}, 16:43, 1977.

\end{thebibliography}

\end{document}